# Relaxation of Non-Convex Integral Functionals in the Multidimensional Scalar Case

Tommaso Bertin, Paulin Huguet

**Abstract**

We study integral functionals defined on scalar Sobolev spaces of the form

$$E[f] : u \mapsto \int_\Omega f(x, u(x), \nabla u(x))dx,$$

with an emphasis on the non-convex case, and the difficulties it involves to prevent the Lavrentiev phenomenon. We determine a formulation of the lower semicontinuous envelope of $E[f]$ with respect to various topologies and with fixed Lipschitz Dirichlet boundary conditions.

# Contents

# Introduction

We consider the following integral functional

$$E[f](u) = \int_\Omega f(x, u(x), \nabla u(x))dx,$$

where $\Omega$ is an open, bounded, Lipschitz subset of $\mathbb{R}^N$, $u \in W^{1,p}(\Omega, \mathbb{R})$ with $1 \leq p \leq +\infty$ and $f : \Omega \times \mathbb{R} \times \mathbb{R}^N \to [0, \infty)$ is a Lagrangian which satisfies suitable hypotheses. Given $\varphi \in W^{1,\infty}(\Omega)$ and $1 \leq p < q \leq +\infty$, we say that the Lavrentiev Phenomenon occurs between $W^{1,p}_\varphi(\Omega)$ and $W^{1,q}_\varphi(\Omega)$ if

$$\inf_{W^{1,p}_\varphi(\Omega)} E[f] < \inf_{W^{1,q}_\varphi(\Omega)} E[f],$$

where $W^{1,p}_\varphi(\Omega)$ denotes the space of $W^{1,p}$ functions on $\Omega$ agreeing with $\varphi$ on $\partial\Omega$. Lavrentiev [37] and Manià [39] proposed the first examples of Lavrentiev Phenomenon. In particular, the example of Manià is a polynomial Lagrangian in the one-dimensional scalar case. Further examples, in particular in the non autonomous case, can be found in [22], [44], [31], [32], [2], [3], [11] and [7].

The study of the non-occurrence of the Lavrentiev Phenomenon, in particular for $p = 1$, $q = \infty$, is important for the application of several numerical approximation techniques, for example the finite elements method. More generally, given $u \in W^{1,p}(\Omega)$, we seek a more regular sequence $(u_n) \subset W^{1,q}(\Omega)$ that converges in some sense to $u$ and such that

$$E[f](u_n) \to E[f](u).$$

This problem has been studied in the one-dimensional case with weak assumptions on the Lagrangian ([1],[41]); usually in the multidimensional scalar case the Lagrangian is assumed to be convex with respect to the gradient variable ([18], [19], [9], [8], [10], [7]). The relation between relaxation of functionals and the Lavrentiev Phenomenon was proposed, as far as we know, for the first time in [17]. This problem has been well studied in the literature, in particular when $f$ is convex with respect to the last variable ([30],[9], [8], [10], [11], [25], [26]). Roughly speaking, the weak sequential lower semicontinuity of $E[f]$ in $W^{1,1}(\Omega)$ is equivalent to the convexity of the Lagrangian with respect to the last variable. If the Lagrangian is no more convex with respect to the gradient variable, it is interesting to study the weak-$*$ lower semicontinuous envelope of $E[f]$ in $W^{1,\infty}_\varphi(\Omega)$. Assuming some continuity assumptions with respect to $(u, \xi)$, in [27] and [40] the following integral representation formula holds: for every $u \in W^{1,\infty}(\Omega)$,

$$\inf\left\{\liminf_{n\to+\infty} \int_\Omega f(x, u_n(x), \nabla u_n(x))dx \,\bigg|\, (u_n) \subset W^{1,\infty}_u(\Omega),\ u_n \rightharpoonup^* u \ \text{ in } \ W^{1,\infty}(\Omega)\right\}$$

$$= \int_\Omega f^{**}(x, u(x), \nabla u(x))dx, \quad (0.1)$$

where $f^{**}$ is the bipolar of the Lagrangian with respect to the last variable.

The integral representation of the lower semicontinuous envelope was studied by many authors, we cite for example [13], [15], [14], [16], [38], [21], [12], [34],[23], starting from [42], [27] and [40]. In [6] the formula (0.1) is proved without any continuity assumptions with respect to the last variable, but under some conditions on the state variable $u$. Furthermore, given $u \in W^{1,\infty}(\Omega)$, there exists $(u_n) \subset u + W_0^{1,\infty}(\Omega)$ such that

$$\|u_n - u\|_{L^\infty} \to 0, \quad \text{and} \quad E[f](u_n) \to E[f^{**}](u). \tag{0.2}$$

At this point, the natural question is whether, given $u \in W^{1,\infty}(\Omega)$, there exists $(u_n) \subset u + W_0^{1,\infty}$ such that

$$u_n \rightharpoonup^* u, \quad \text{and} \quad E[f](u_n) \to E[f^{**}](u), \tag{0.3}$$

which would be stronger than (0.2). To be more precise, along all this article, we will use that that a sequence $(u_n) \subset W^{1,\infty}(\Omega)$ goes weakly-$*$ to some $u$, denoted by $u_n \rightharpoonup^* u$ if

$$\|u_n - u\|_{L^\infty} \to 0 \quad \text{and} \quad \sup_n \|\nabla u_n\|_{L^\infty} < +\infty. \tag{0.4}$$

In Section 1, we expose a condition that implies the existence of an approximating sequence as in (0.3). The main novelty compared to other similar results, such as in [27] or [40], is that the assumed regularity of the Lagrangian is weaker than usual (especially with respect to the last variable) and that the convergence considered here is with respect to the weak-$*$ topology in $W^{1,\infty}(\Omega)$ (compared to uniform convergence in already existing literature). It is referenced as Theorem 1.22 in the following and is stated as such:

**Theorem 1.22.** *Let $f : \Omega \times \mathbb{R} \times \mathbb{R}^N \to [0, \infty)$ satisfy Hypothesis $\mathcal{A}$. Assume that $f$ also satisfies condition $(\mathcal{K})$. Then for every $u \in W^{1,\infty}(\Omega)$, there exists a sequence $(u_n) \subset W_u^{1,\infty}(\Omega)$ such that*
$$u_n \rightharpoonup^* u \quad in \quad W^{1,\infty}(\Omega),$$
*and*
$$\lim_{n \to +\infty} \int_\Omega f(x, u_n(x), \nabla u_n(x))dx = \int_\Omega f^{**}(x, u(x), \nabla u(x))dx.$$

We do not detail conditions $\mathcal{A}$ and $(\mathcal{K})$ here. Let us only mention that Hypothesis $\mathcal{A}$ refers to regularity assumptions on $f$, while condition $(\mathcal{K})$ is a geometrical condition introduced in Section 1.3. The main appeal of this new result compared to the already existing literature is condition $(\mathcal{K})$, which as stated before, ensure the existence of a converging sequence in $W^{1,\infty}(\Omega)$ weak-$*$. We find in Corollaries 1.26 and 1.30 two distinct conditions which are sufficient for this condition $(\mathcal{K})$ to hold. The first condition is the uniform boundedness of the connected components of the detachment set (i.e., the set where $f^{**} < f$); the second condition is the superlinearity in the variable $\xi$ uniformly in the other two variables.

The main idea of the proofs is to show that, under one of these conditions, one has

$$\inf\left\{\liminf_{n\to+\infty}\int_\Omega f(x,u_n(x),\nabla u_n(x))dx \,\Big|\, u_n \rightharpoonup^* u \quad \text{in} \quad W^{1,\infty}(\Omega)\right\}$$
$$=\inf\left\{\liminf_{n\to+\infty}\int_\Omega f(x,u_n(x),\nabla u_n(x))dx \,\Big|\, u_n \rightharpoonup^* u \quad \text{in} \quad W^{1,\infty}(\Omega), \quad \|\nabla u_n\|_{L^\infty} < K\right\},$$

where $u_n = u$ on $\partial\Omega$ in the trace sense and $K$ is chosen large enough so that $\|\nabla u\|_{L^\infty} < K$.

In Section 2, we focus on the autonomous case (that is, no dependency of the Lagrangian on the $x$ variable). In particular, we apply the existence of an approximating sequence in $W^{1,\infty}(\Omega)$ that satisfies (0.2) to the nonconvex case. Firstly, we use Theorem 1.19 to avoid the Lavrentiev phenomenon for a large class of Lagrangian (see Theorem 2.3). Then we prove the following

**Theorem 2.4.** *Let $f : \mathbb{R} \times \mathbb{R}^N \to [0,\infty)$ satisfy Hypothesis B, and let $\varphi \in W^{1,\infty}(\Omega)$. Then for any $u \in W^{1,1}_\varphi(\Omega)$, there exists a sequence $(u_n) \subset W^{1,\infty}_\varphi(\Omega)$ such that*

$$u_n \to u \quad \text{in} \quad L^1(\Omega),$$

*and*

$$\lim_{n\to+\infty}\int_\Omega f(u_n(x),\nabla u_n(x))dx = \int_\Omega f^{**}(u(x),\nabla u(x))dx.$$

*Moreover,*

- *if there exists $\Phi : \mathbb{R}^N \to [0,\infty)$ superlinear such that*

$$f(u,\xi) \geq \Phi(\xi), \qquad \forall u \in \mathbb{R}, \forall \xi \in \mathbb{R}^N,$$

  *then the sequence $(u_n)$ can be chosen so that $u_n \rightharpoonup u$ weakly in $W^{1,1}(\Omega)$;*

- *if for some $p \in (1,\infty)$ it holds that*

$$f(u,\xi) \geq c_1|\xi|^p - c_2, \qquad \forall u \in \mathbb{R}, \forall \xi \in \mathbb{R}^N,$$

  *for some $c_1, c_2 > 0$, and $u \in W^{1,p}(\Omega)$ then the sequence $(u_n)$ can be taken so that $u_n \rightharpoonup u$ weakly in $W^{1,p}(\Omega)$.*

The main idea is to apply the result [8, Theorem 1.1] to $E[f^{**}]$ to find an approximating sequence $(u_k) \subset W^{1,\infty}_\varphi(\Omega)$ such that $u_k \to u$ in $W^{1,1}(\Omega)$ and $E[f^{**}](u_k) \to E[f^{**}](u)$. For every $u_k$, we can construct an approximating sequence as in (0.2), then use a diagonal extraction argument. Furthermore, in Theorem 2.8, if $E[f](u) = E[f^{**}](u)$ we recover the strong convergence in $W^{1,p}(\Omega)$ of the approximating sequence:

**Theorem 2.8.** *Let $f : \mathbb{R} \times \mathbb{R}^N \to [0, \infty)$ satisfies Hypothesis $\mathcal{B}$, as well as $\varphi \in W^{1,\infty}(\Omega)$. Assume $u \in W^{1,p}_\varphi(\Omega)$ for some $p \in [1, \infty)$, and satisfies*

$$\int_\Omega f^{**}(u(x), \nabla u(x))dx = \int_\Omega f(u(x), \nabla u(x))dx.$$

*Then there exists a sequence $(u_n) \subset W^{1,\infty}_\varphi(\Omega)$ such that*

$$u_n \to u \quad \text{strongly in} \quad W^{1,p}(\Omega),$$

*and*

$$\lim_{n \to +\infty} \int_\Omega f(u_n(x), \nabla u_n(x))dx = \int_\Omega f(u(x), \nabla u(x))dx.$$

We remark that, in all the results presented above, the Lagrangian is not necessarily assumed to be continuous with respect to the variable $\xi$. If the Lagrangian is continuous and is dominated by a convex function $g$ such that $E[g](u) < +\infty$, (Theorem 2.13 and corollaries) we can find a sequence $(u_n) \subset W^{1,\infty}(\Omega)$ such that $u_n \to u$ strongly in $W^{1,p}(\Omega)$ and $E[f](u_n) \to E[f](u)$. In fact, all of these results also ensure the conservation of the boundary condition all along the approximating sequence, assuming that $u_{|\partial\Omega} \in \text{Lip}(\partial\Omega)$. The main idea is to use [8, Theorem 1.1] to the dominating function and use Fatou Lemma. We can also see these results as integral representations of the lower semicontinuous envelope with respect to the strong topology of $W^{1,p}(\Omega)$, this is the content of Remark 2.6.

Then in Section 3, we apply a result in [7] to extend the results of Section 2 to the non-autonomous case, assuming an anti-jump condition with respect to the variable $x$ (called condition $(\mathcal{H}_1)$ in this paper). The main novelty is to prove that condition $(\mathcal{H}_1)$ for a non convex Lagrangian implies that the condition holds also for its bipolar. At this point, most of the proofs follow similar patterns to the ones in Section 2. The main additional difficulty is to make sure that condition $(\mathcal{H}_1)$ stays true for the various auxiliary Lagrangians we consider for calculation purpose.

# 1 Geometrical conditions for the weak-$*$ approximation in $W^{1,\infty}(\Omega)$

This section is devoted to finding geometric conditions on the Lagrangian $f$ such that for every $u \in W^{1,\infty}(\Omega)$ there exists a sequence $(u_n) \subset u + W^{1,\infty}_0(\Omega)$ such that

$$u_n \rightharpoonup^* u \quad \text{in} \quad W^{1,\infty}(\Omega), \tag{1.1}$$

and

$$\int_\Omega f(x, u_n(x), \nabla u_n(x))dx \to \int_\Omega f^{**}(x, u(x), \nabla u(x))dx. \tag{1.2}$$

Our main result and contribution for this section will be Theorem 1.22, stating that, under a geometric condition on the Lagrangian (see condition $(\mathcal{K})$), a sequence such that (1.1) and (1.2) hold does exist. We will then, in Section 1.4 work to state a few sufficient properties on the Lagrangian $f$ ensuring that it satisfy $(\mathcal{K})$.

**Notation 1.1.** In the whole article, we will use the following conventions and notations:

- $\Omega$ is an open, bounded, Lipschitz domain of $\mathbb{R}^N$, for some $N \geq 1$.

- Given $K > 0$, $B_K$ will denote the closed ball of radius $K$ and center $0$ in $\mathbb{R}^N$.

- Given $\varphi \in W^{1,\infty}(\Omega)$, we will write

$$W^{1,p}_\varphi(\Omega) = \varphi + W^{1,p}_0(\Omega) = \Big\{ u \in W^{1,p}(\Omega) \ : \ u = \varphi \quad \text{on} \quad \partial\Omega \Big\},$$

  for every $p \in [1, \infty]$. Hence, $W^{1,p}_\varphi(\Omega)$ is the set of functions in $W^{1,p}(\Omega)$ which agree with $\varphi$ (in the sense of the trace) on $\partial\Omega$.

## 1.1 On measurability and basic assumptions

**Definition 1.2.** Given a metric space $X$, we denote by $\mathcal{B}(X)$ the Borel $\sigma$-algebra of $X$. Moreover, if $X$ is a subset of an euclidean space, then $\mathcal{L}(X)$ will denote its Lebesgue $\sigma$-algebra. We will say that a function $f : \Omega \times \mathbb{R} \times \mathbb{R}^N \to [0, \infty]$ is Lebesgue-Borel measurable if it is measurable for the $\sigma$-algebras

$$\mathcal{L}(\Omega) \otimes \mathcal{B}(\mathbb{R}) \otimes \mathcal{B}(\mathbb{R}^N) \longrightarrow \mathcal{B}([0, \infty]).$$

This notion of measurability for Lagrangians is made in such a way that, if $u : \Omega \to \mathbb{R}$ and $v : \Omega \to \mathbb{R}^N$ are Lebesgue measurable, then

$$x \mapsto f(x, u(x), v(x)),$$

is measurable on $\Omega$, and thus the quantity $\int_\Omega f(x, u(x), v(x))dx$ makes sense. This will of course be used in the case $u \in W^{1,1}(\Omega)$ and $v = \nabla u$. In a large part of the literature on the subject, the usual assumptions made on the Lagrangian is the following Carathéodory property:

**Definition 1.3.** We say that $f : \Omega \times \mathbb{R} \times \mathbb{R}^N \to [0, \infty]$ is a Carathéodory function if it is Lebesgue-Borel measurable and for a.e. $x \in \Omega$, the mapping $(u, \xi) \mapsto f(x, u, \xi)$ is continuous on $\mathbb{R}^N$.

In this paper, we will usually have a less restrictive assumption on the Lagrangian. In particular, the continuity in $\xi$ will be withdrawn from the assumptions for most of Section 1. We assume the following for the Lagrangian $f$:

**Hypothesis $\mathcal{A}$.** The function $f : \Omega \times \mathbb{R} \times \mathbb{R}^N \to [0, \infty)$ satisfies

a) $f$ is Lebesgue-Borel measurable;

b) for a.e. $x \in \Omega$, the function $u \mapsto f(x, u, \xi)$ is continuous with respect to $u$ uniformly as $\xi$ varies in bounded sets. That is, for every bounded set $B \subset \mathbb{R}^N$, for every $u_0 \in \mathbb{R}$,

$$\forall \varepsilon > 0, \ \exists \delta > 0, \ \forall u \in \mathbb{R}, \forall \xi \in B,$$
$$|u - u_0| < \delta \quad \Rightarrow \quad |f(x, u, \xi) - f(x, u_0, \xi)| < \varepsilon;$$

c) for every bounded set $B \subset \mathbb{R} \times \mathbb{R}^N$, there exists $a \in L^1(\Omega)$ such that $f(x, u, \xi) \leq a(x)$ for a.e. $x \in \Omega$ and all $(u, \xi) \in B$;

d) for every $u \in W^{1,\infty}(\Omega)$, for every bounded set $B \subset \mathbb{R}^N$ and for every $\eta > 0$ there exists $T \subset \Omega$ compact such that $|\Omega \setminus T| < \eta$ and $x \mapsto f(x, u(x), \xi)$ is continuous on $T$ uniformly as $\xi$ varies in $B$, that is, for every $x_0 \in T$,

$$\forall \varepsilon > 0, \ \exists \delta > 0, \ \forall x \in T, \forall \xi \in B,$$
$$|x - x_0| < \delta \quad \Rightarrow \quad |f(x, u(x), \xi) - f(x_0, u(x_0), \xi)| < \varepsilon.$$

**Remark 1.4.** We point out the fact that we can replace Hypotheses $\mathcal{A}$-b) and -d) with the more restrictive request that for every bounded set $B \subset \mathbb{R}^N$ and for every $\eta > 0$ there exists a compact set $T \subset \Omega$ such that $|\Omega \setminus T| < \eta$ and $f$ is continuous with respect to $(x, u) \in T \times \mathbb{R}$ uniformly as $\xi$ varies in $B$.

A typical example of a Lagrangian $f$ satisfying Hypothesis $\mathcal{A}$, is if the dependency in $\xi$ is bounded and separated from the other variables. That is, if

$$f(x, u, \xi) = g(x, u)h(\xi),$$

for some $g : \Omega \times \mathbb{R} \to [0, \infty)$ Carathéodory satisfying that for every compact interval $I \subset \mathbb{R}$, there exists $a \in L^1(\Omega)$ such that $g(x, u) \leq a(x)$, for a.e. $x \in \Omega$ and all $u \in I$; and $h : \mathbb{R}^N \to [0, \infty)$ Borel bounded on bounded set. Another interesting case, which is very useful in application is for Carathéodory Lagrangians:

**Proposition 1.5.** *Let* $f : \Omega \times \mathbb{R} \times \mathbb{R}^N \to [0, \infty)$ *be a Carathéodory function. Then* $f$ *satisfies assumptions* $\mathcal{A}$*-a), -b) and -d).*

*Proof.* $\mathcal{A}$-a) is obvious by definition. To prove $\mathcal{A}$-b), it is enough to use the fact that, for a.e. $x \in \Omega$, the map $(u, \xi) \mapsto f(x, u, \xi)$ is continuous. The rest follows from a simple compactness argument on the bounded set $B$.

We now prove $\mathcal{A}$-d). Let $u \in W^{1,\infty}(\Omega)$, $B \subset \mathbb{R}^N$ a bounded set and $\eta > 0$. Let $g(x, \xi) := f(x, u(x), \xi)$. Then $g$ is $\mathcal{L}(\Omega) \otimes \mathcal{B}(\mathbb{R}^N)$-measurable and continuous with respect to $\xi$. By the Scorza-Dragoni Theorem (see [27, Chapter VIII, Section 1.3]), there exists $T \subset \Omega$ compact such that $|\Omega \setminus T| < \eta$ and $g$ is continuous on $T \times \mathbb{R}^N$. Now since $B$ is bounded, another compactness argument shows that $g$ is continuous with respect to $x \in T$ uniformly as $\xi$ varies in $B$. □

**Remark 1.6.** In this article, we will work only with non-negative Lagrangians for simplicity. However, it might be interesting to keep in mind that the various approximation results in Section 1.3 (Proposition 1.17, and Theorems 1.19, 1.22) would still hold if $f$ were only assumed to be real-valued, and with assumption $\mathcal{A}$-c) replaced by c') for every bounded set $B \subset \mathbb{R} \times \mathbb{R}^N$, there exists $a \in L^1(\Omega)$ such that $|f(x, u, \xi)| \leq a(x)$ for a.e. $x \in \Omega$ and all $(u, \xi) \in B$.

Indeed, if c') were to hold, considering the (non-negative) auxiliary Lagrangian $g(x, u, \xi) = f(x, u, \xi) + a(x)$, it would satisfy Hypothesis $\mathcal{A}$. All of these approximation results would thus hold for $g$, and thus for $f$ as an immediate calculation would show.

## 1.2 Basic tools, notations and some reminder about convexification

We introduce the following technical definition:

**Definition 1.7.** For every $f : \Omega \times \mathbb{R} \times \mathbb{R}^N \to [0, \infty]$ and for every $T \subset \mathbb{R}^N$ we define

$$\tilde{f}_T(x, u, \xi) := \begin{cases} f(x, u, \xi) & \text{if } \xi \in T \\ +\infty & \text{otherwise.} \end{cases} \quad (1.3)$$

In the special case $T = B_K$ for some $K > 0$, we let

$$\tilde{f}_K := \tilde{f}_{B_K}.$$

As convexification is a crucial tool whenever mentioning weak or weak-$*$ convergence for integral functional, we define it and mention some standard facts on the subject.

**Definition 1.8.** Let $f : \mathbb{R}^N \to [0, \infty]$. We define the function $f^{**} : \mathbb{R}^N \to [0, \infty]$ as the convexification of $f$. That is, $f^{**}$ is the greatest lower semicontinuous convex function on $\mathbb{R}^N$ which stays under $f$.

The following facts are classical, the reader can have a look at [27, Chapter I] for a more complete overview of convex analysis.

**Proposition 1.9.** *Let $f : \mathbb{R}^N \to [0, \infty]$.*

- *The function $f^{**}$ can be written as the supremum of all affine maps on $\mathbb{R}^N$ which bound $f$ from below.*

- *Assume that $f$ takes finite values in a neighborhood of some $\xi \in \mathbb{R}^N$. Then*

$$f^{**}(\xi) = \inf\left\{\sum_i \alpha_i f(\xi_i)\right\}, \quad (1.4)$$

  *where the infimum is taken over all convex combinations $(\alpha_i, \xi_i) \subset [0, 1] \times \mathbb{R}^N$ such that $\sum_i \alpha_i = 1$ and $\xi = \sum_i \alpha_i \xi_i$.*

*Proof.* The first point is a consequence of the Hahn-Banach theorem (see [27, Chapter I, Proposition 3.1] for a detailed proof). We prove the second point. Let $g : \mathbb{R}^N \to [0, \infty]$ be the functional on the right hand side of (1.4). We claim that $g$ is convex. Indeed, take $\xi, \zeta \in \mathbb{R}^N$ as well as $\lambda \in (0, 1)$. Let $\sum_i \alpha_i \xi_i$ and $\sum_j \beta_j \zeta_j$ be convex combinations of $\xi$ and $\zeta$ respectively. Then $\sum_i (1-\lambda)\alpha_i \xi_i + \sum_j \lambda \beta_j \zeta_j$ is a convex combination of $(1-\lambda)\xi + \lambda\zeta$, and by definition of $g$,

$$g((1-\lambda)\xi + \lambda\zeta) \leq (1-\lambda)\sum_i \alpha_i f(\xi_i) + \lambda \sum_j \beta_j f(\zeta_j).$$

Now taking the infimum over all convex combinations of $\xi$ and $\zeta$ respectively, we finally obtain $g((1-\lambda)\xi + \lambda\zeta) \leq (1-\lambda)g(\xi) + \lambda g(\zeta)$. Now using the fact that $f^{**}$ is convex and $f^{**} \leq f$, for every $\xi$,

$$f^{**}(\xi) = \inf\left\{\sum_i \alpha_i f^{**}(\xi_i)\right\} \leq \inf\left\{\sum_i \alpha_i f(\xi_i)\right\} = g(\xi),$$

where the infimum is again taken over all convex combinations of $\xi$. Therefore $f^{**} \leq g \leq f$. Now assume that $f$ is finite in a neighborhood of some $\xi \in \mathbb{R}^N$. By [27, Chapter I, Corollary 2.3], $g$ is continuous at point $\xi$ and by [27, Chapter I, Proposition 5.2], there exists an affine map $\ell$ on $\mathbb{R}^N$ such that

$$\ell \leq g \leq f, \quad \text{and} \quad \ell(\xi) = g(\xi).$$

By the first point of the proposition, one has $f^{**}(\xi) \geq \ell(\xi) = g(\xi)$, which achieves the proof. □

**Remark 1.10.** Using the notation introduced in (1.3), $(\tilde{f}_T)^{**}$ is the greatest function which is convex, lower semicontinuous with respect to the last variable and stays under $f$ on $T$. In particular, $(\tilde{f}_K)^{**} \geq f^{**}$ on $B_K$ and in general the strict inequality may hold. Indeed, consider for instance $f(\xi) = \exp(-|\xi|)$. Then $f^{**} \equiv 0$ on $\mathbb{R}^N$, but for every $K > 0$,

$$(\tilde{f}_K)^{**} = e^{-K} \quad \text{on} \quad B_K.$$

The previous remark allows to formulate the following lemma:

**Lemma 1.11.** *Let $f : \mathbb{R}^N \to [0,\infty)$. Then the family of mapping $((\tilde{f}_K)^{**})_{K>0}$ is non-increasing as $K \to +\infty$ and converge pointwise to $f^{**}$.*

*Proof.* The fact that the family is non-increasing as $K \to +\infty$ easily follows from the definition of $(\tilde{f}_K)^{**}$. Now define $g := \lim_{K \to +\infty}(\tilde{f}_K)^{**} = \inf_{K>0}(\tilde{f}_K)^{**}$. By definition of $f^{**}$, it is clear that $f^{**} \leq g$ (see Remark 1.10). Moreover, since the family $((\tilde{f}_K)^{**})_{K>0}$ is totally ordered, it holds that $g$ is convex and $g \leq f$ on all of $\mathbb{R}^N$ (given $\xi \in \mathbb{R}^N$, one can choose $K > |\xi|$). Since $f$ is finite valued, so is $g$ and thus $g$ is lower semicontinuous on $\mathbb{R}^N$. Therefore $g \leq f^{**}$. □

In the case of an integral functional which may depend on $x$ and $u$ and not just on $\nabla u$, it is relevant to consider the convexification only with respect to the variable $\xi$:

**Definition 1.12.** Let $f : \Omega \times \mathbb{R} \times \mathbb{R}^N \to [0,\infty]$. We define the function $f^{**} : \Omega \times \mathbb{R} \times \mathbb{R}^N \to [0,\infty]$ as the convexification of $f$ with respect to the last variable. More specifically, for every $(x,u) \in \Omega \times \mathbb{R}$, $f^{**}(x,u,\cdot)$ is the greatest lower semicontinuous convex function on $\mathbb{R}^N$ which stays under $f(x,u,\cdot)$.

Since the Lagrangian $f$ might be non-convex with respect to the last variable, we prove that in the case where $f$ is finite-valued (which is the only relevant case for us), $f^{**}$ is at least measurable.

**Lemma 1.13.** *Let $f : \Omega \times \mathbb{R} \times \mathbb{R}^N \to [0, \infty)$ satisfy Hypothesis $\mathcal{A}$. Then for every $K > 0$, $(\tilde{f}_K)^{**}$ is a Carathéodory function on $\Omega \times \mathbb{R} \times (B_K)^\circ$. Furthermore $f^{**}$ is Lebesgue-Borel measurable on $\Omega \times \mathbb{R} \times \mathbb{R}^N$.*

Here $(B_K)^\circ$ denote the centered open ball of radius $K$.

*Proof.* Firstly, we prove that $u \mapsto (\tilde{f}_K)^{**}(x, u, \xi)$ is continuous uniformly as $\xi$ varies in $B_K$, for a.e. $x \in \Omega$ and for every $K > 0$. We fix $u_0 \in \mathbb{R}$ and Hypothesis $\mathcal{A}$-b) guarantees that for a.e. $x \in \Omega$, for every $\varepsilon > 0$, there exists $\delta > 0$ such that if $|u - u_0| < \delta$ then

$$|f(x, u, \xi) - f(x, u_0, \xi)| < \varepsilon, \qquad \forall \xi \in B_K, \tag{1.5}$$

and so we have

$$(\tilde{f}_K)^{**}(x, u, \xi) - \varepsilon \leq f(x, u, \xi) - \varepsilon \leq f(x, u_0, \xi), \qquad \forall \xi \in B_K. \tag{1.6}$$

Now $(\tilde{f}_K)^{**}(x, u, \xi) - \varepsilon$ is convex and lower semicontinuous in $\xi$ and so

$$(\tilde{f}_K)^{**}(x, u, \xi) - \varepsilon \leq (\tilde{f}_K)^{**}(x, u_0, \xi), \qquad \forall \xi \in B_K.$$

Reversing the roles of $u$ and $u_0$ in (1.6), one gets in the end

$$|(\tilde{f}_K)^{**}(x, u, \xi) - (\tilde{f}_K)^{**}(x, u_0, \xi)| \leq \varepsilon, \qquad \forall \xi \in B_K. \tag{1.7}$$

Now for every $u \in \mathbb{R}$, since the function $f(x, u, \cdot)$ is finite-valued, it holds that (see [27, Chapter I, Corollary 2.3])

$$(\tilde{f}_K)^{**}(x, u, \cdot) \quad \text{is continuous on} \quad (B_K)^\circ. \tag{1.8}$$

Therefore, if $(u_n, \xi_n) \subset \mathbb{R} \times (B_K)^\circ$ is a sequence converging to $(u_0, \xi_0) \in \mathbb{R} \times (B_K)^\circ$, then

$$\limsup_{n \to +\infty} |(\tilde{f}_K)^{**}(x, u_n, \xi_n) - (\tilde{f}_K)^{**}(x, u_0, \xi_0)|$$
$$\leq \limsup_{n \to +\infty} |(\tilde{f}_K)^{**}(x, u_n, \xi_n) - (\tilde{f}_K)^{**}(x, u_0, \xi_n)|$$
$$+ \limsup_{n \to +\infty} |(\tilde{f}_K)^{**}(x, u_0, \xi_n) - (\tilde{f}_K)^{**}(x, u_0, \xi_0)|$$
$$= 0,$$

by (1.7) and (1.8). That is, for a.e. $x \in \Omega$, $(\tilde{f}_K)^{**}(x, \cdot, \cdot)$ is continuous in $\mathbb{R} \times (B_K)^\circ$. By [27, Chapter VIII, Proposition 1.1], to prove that $(\tilde{f}_K)^{**}$ is Carathéodory, it is therefore enough to show that $x \mapsto (\tilde{f}_K)^{**}(x, u, \xi)$ is Lebesgue measurable for every $(u, \xi) \in \mathbb{R} \times (B_K)^\circ$. We fix $(u_0, \xi_0) \in \mathbb{R} \times (B_K)^\circ$. Using $\mathcal{A}$-d), there exists an increasing sequence of compact sets $(T_n)$ of $\Omega$ such that $|\Omega \setminus T_n| \to 0$ and for every $n$, the mapping $x \mapsto f(x, u_0, \xi)$ is continuous on $T_n$ uniformly as $\xi$ varies in $B_K$. Now by a similar argument to the one developed above (equations (1.5)-(1.8)), $x \mapsto (\tilde{f}_K)^{**}(x, u_0, \xi_0)$ is continuous on $T_n$. Taking the limit as $n \to +\infty$, this same mapping is an almost everywhere limit of continuous maps, and is thus Lebesgue measurable on $\Omega$.

Now, for every $K \in \mathbb{N}$, we let $g_K$ be defined as such:

$$g_K(x, u, \xi) = \begin{cases} (\tilde{f}_K)^{**}(x, u, \xi) & \text{if } |\xi| \leq K/2 \\ +\infty & \text{if } |\xi| > K/2. \end{cases}$$

Since $(\tilde{f}_K)^{**}$ is Carathéodory on $\Omega \times \mathbb{R} \times (B_K)^\circ$, then $g_K$ is Lebesgue-Borel measurable on $\Omega \times \mathbb{R} \times \mathbb{R}^N$. Recalling Lemma 1.11, we have

$$f^{**}(x, u, \xi) = \inf_K g_K(x, u, \xi) = \lim_{K \to +\infty} g_K(x, u, \xi),$$

hence, $f^{**}$ is Lebesgue-Borel measurable on $\Omega \times \mathbb{R} \times \mathbb{R}^N$. □

**Remark 1.14.** It might be interesting to point out that $f^{**}$ may not be a Carathéodory function in full generality. The issue is that the joint continuity in $(u, \xi)$ may not be satisfied. We give here an example, which was initially presented in [40, Example 3.11]. Consider the following mapping $f(u, \xi) := (|\xi| + 1)^{|u|}$ for $(u, \xi) \in \mathbb{R} \times \mathbb{R}^1$. It is clear that it is continuous, however

$$f^{**}(u, \xi) = \begin{cases} f(u, \xi) & \text{if } |u| \geq 1 \\ 1 & \text{if } |u| < 1. \end{cases}$$

which is not continuous (not even lower semicontinuous).

However, we will now study some conditions on the Lagrangian $f$ to ensure that $f^{**}$ is Carathéodory. The problem has already been studied in [40, Corollary 3.12]. Assuming a $p$-growth or continuity in $u$ uniformly in $\xi$, the authors proved that $f^{**}$ is a Carathéodory function. Theorem 1.16 below, with the help of a different kind of argument, extends the $p$-growth case to superlinear Lagrangians.

**Proposition 1.15.** *Let $f : \mathbb{R} \times \mathbb{R}^N \to [0, \infty)$ satisfies Hypothesis A. If for every compact interval $I \subset \mathbb{R}$ and $K > 0$, there exists $K' \geq K$ such that*

$$(\tilde{f}_{K'})^{**} = f^{**} \quad \text{on} \quad I \times B_K, \tag{1.9}$$

*then $f^{**}$ is continuous on $\mathbb{R} \times \mathbb{R}^N$.*

*Proof.* It is enough to prove that $f^{**}$ is continuous on $I \times B_K$ for every compact interval $I \subset \mathbb{R}$ and $K > 0$. Fix such a choice of $I$ and $K$. Let $K' \geq K$ such that (1.9) apply. If $K' > K$, by Lemma 1.13, $(\tilde{f}_{K'})^{**}$ is continuous on $I \times B_K$, and (1.9) gives our conclusion. If $K = K'$, we cannot apply immediately Lemma 1.13 because $(\tilde{f}_{K'})^{**}$ may not be continuous on $\{\xi \in \partial B_K\}$. However, choosing $K'' = K' + 1 > K$, then

$$(\tilde{f}_{K'})^{**} \geq (\tilde{f}_{K''})^{**} \geq f^{**} = (\tilde{f}_{K'})^{**} \quad \text{on} \quad I \times B_K.$$

Thus $(\tilde{f}_{K''})^{**} = f^{**}$ on $I \times B_K$, and since $(\tilde{f}_{K''})^{**}$ is continuous on $I \times B_K$ by Lemma 1.13, this proves the result in the case $K' = K$. □

The condition (1.9) used in the statement of the previous proposition is quite similar to condition ($\mathcal{K}$) presented below in Section 1.3, which will play a major role in our discussion. Actually, a similar argument given in Remark 1.21 was used in the proof of Proposition 1.15.

We also introduce the following result, which ensure that $f^{**}$ is Carathéodory.

**Theorem 1.16.** *Let* $f : \Omega \times \mathbb{R} \times \mathbb{R}^N \to [0, \infty)$ *satisfy Hypothesis $\mathcal{A}$. Assume that for a.e.* $x \in \Omega$ *and every compact interval* $I \subset \mathbb{R}$*, there exists* $\Phi_{x,I} : \mathbb{R}^N \to [0, \infty)$ *superlinear such that*

$$f(x, u, \xi) \geq \Phi_{x,I}(\xi), \qquad \forall u \in I, \ \forall \xi \in \mathbb{R}^N. \tag{1.10}$$

*Then $f^{**}$ is a Carathéodory function.*

Before proving it rigorously, we will need some other tools, in particular regarding superlinearity (Theorem 1.27). We refer the reader to the end of Section 1.4 for a detailed proof of Theorem 1.16.

## 1.3 Known results of weak-$*$ relaxation on $W^{1,\infty}(\Omega)$ and condition ($\mathcal{K}$).

We now present multiple results of relaxation on $W^{1,\infty}(\Omega)$ with respect to the weak-$*$ topology. As mentioned in the introduction, the main contribution of this section is Theorem 1.22. Before that, we need some intermediary results, already proved (for instance by P. Marcellini and C. Sbordone in [40]). The following proposition, proved in [6], is the most recent formulation, with the weakest initial assumptions that we are aware of.

**Proposition 1.17.** *Let* $f : \Omega \times \mathbb{R} \times \mathbb{R}^N \to [0, \infty)$ *satisfy Hypothesis $\mathcal{A}$. For every* $u \in W^{1,\infty}(\Omega)$ *and for every* $K \in \mathbb{N}$ *such that*

$$\|\nabla u\|_{L^\infty} < K, \tag{1.11}$$

*there exists a sequence* $(u_{K,n}) \subset W_u^{1,\infty}(\Omega)$ *such that*

$$\|\nabla u_{K,n}\|_{L^\infty} < K, \qquad u_{K,n} \rightharpoonup^* u \quad in \quad W^{1,\infty}(\Omega),$$

*and*

$$\lim_{n \to +\infty} \int_\Omega f(x, u_{K,n}(x), \nabla u_{K,n}(x)) dx = \int_\Omega (\tilde{f}_K)^{**}(x, u(x), \nabla u(x)) dx.$$

Here, the weak-$*$ convergence denoted by $\rightharpoonup^*$ refers to the one introduced in (0.4). We recall that, consistent with Notation 1.1, $W_u^{1,\infty}(\Omega)$ denotes the space of functions in $W^{1,\infty}(\Omega)$ which agree with $u$ on $\partial\Omega$.

By the previous result we have that for every $K \in \mathbb{N}$, $u \in W^{1,\infty}(\Omega)$ satisfying (1.11),

$$\inf \left\{ \liminf_{n \to +\infty} \int_\Omega f(x, u_n(x), \nabla u_n(x)) dx \ \middle| \ \begin{array}{c} (u_n) \subset W_u^{1,\infty}(\Omega) \\ u_n \rightharpoonup^* u \\ \|\nabla u_n\|_{L^\infty} < K \end{array} \right\} \leq \int_\Omega (\tilde{f}_K)^{**}(x, u(x), \nabla u(x)) dx. \tag{1.12}$$

Since $(\tilde{f}_K)^{**}$ is Carathéodory and is convex, lower semicontinuous with respect to $\xi$ we have the following (see for instance [36, Chapter 4, Theorem 4.5]):

**Lemma 1.18** (Tonelli). *The integral functional*

$$u \mapsto \int_\Omega (\tilde{f}_K)^{**}(x, u(x), \nabla u(x))dx,$$

*is sequentially lower semicontinuous with respect to the weak topology of $W^{1,1}(\Omega)$.*

So a fortiori, since any sequence converging weakly-$*$ in $W^{1,\infty}(\Omega)$ also converges weakly in $W^{1,1}(\Omega)$, we have that

$$\inf \left\{ \liminf_{n\to+\infty} \int_\Omega f(x, u_n(x), \nabla u_n(x))dx \;\middle|\; \begin{array}{c} (u_n) \subset W^{1,\infty}_u(\Omega) \\ u_n \rightharpoonup^* u \\ \|\nabla u_n\|_{L^\infty} < K \end{array} \right\} \geq \int_\Omega (\tilde{f}_K)^{**}(x, u(x), \nabla u(x))dx. \tag{1.13}$$

(Here we used the fact that $f \geq (\tilde{f}_K)^{**}$ on $\Omega \times \mathbb{R} \times B_K$). Thus, combining (1.12) and (1.13),

$$\min \left\{ \liminf_{n\to+\infty} \int_\Omega f(x, u_n(x), \nabla u_n(x))dx \;\middle|\; \begin{array}{c} (u_n) \subset W^{1,\infty}_u(\Omega) \\ u_n \rightharpoonup^* u \\ \|\nabla u_n\|_{L^\infty} < K \end{array} \right\} = \int_\Omega (\tilde{f}_K)^{**}(x, u(x), \nabla u(x))dx. \tag{1.14}$$

The fact that the infimum is in fact a minimum in (1.14) comes from the existence of a minimizing sequence (see Proposition 1.17). Now $x \mapsto (\tilde{f}_K)^{**}(x, u, \nabla u)$ is a non-increasing sequence bounded from above and converging pointwise to $x \mapsto f^{**}(x, u, \nabla u)$ (by Lemma 1.11). Thus taking the limit as $K \to +\infty$ and using a diagonal argument, we obtain the following result, stated in [40] in the continuous case, or in [6] for the general case:

**Theorem 1.19.** *Let $f : \Omega \times \mathbb{R} \times \mathbb{R}^N \to [0, \infty)$ satisfy Hypothesis A. For every $u \in W^{1,\infty}(\Omega)$, there exists a sequence $(u_n) \subset W^{1,\infty}_u(\Omega)$ such that*

$$u_n \to u \quad in \quad L^\infty(\Omega),$$

*and*

$$\lim_{n\to+\infty} \int_\Omega f(x, u_n(x), \nabla u_n(x))dx = \int_\Omega f^{**}(x, u(x), \nabla u(x))dx.$$

*Furthermore,*

$$\inf \left\{ \liminf_{n\to+\infty} \int_\Omega f(x, u_n(x), \nabla u_n(x))dx \;\middle|\; \begin{array}{c} (u_n) \subset W^{1,\infty}_u(\Omega) \\ u_n \rightharpoonup^* u \end{array} \right\} = \int_\Omega f^{**}(x, u(x), \nabla u(x))dx. \tag{1.15}$$

Notice that (1.15) is obtained by taking the limit as $K \to +\infty$ in (1.14). Now our aim is to find some conditions on the Lagrangian $f$ in order to reach the infimum in (1.15), that is, our question is whether there exists a sequence $(u_n) \subset W^{1,\infty}_u(\Omega)$ such that

$$u_n \rightharpoonup^* u \quad in \quad W^{1,\infty}(\Omega), \tag{1.16}$$

and
$$\int_\Omega f(x, u_n(x), \nabla u_n(x))dx \to \int_\Omega f^{**}(x, u(x), \nabla u(x))dx. \tag{1.17}$$

We will see that in general, such a sequence does not exist. Before moving on, we study the situation in more details. Let $u \in W^{1,\infty}(\Omega)$ and $(u_n) \subset W^{1,\infty}(\Omega)$ be such that $u_n \rightharpoonup^* u$. Let $K' > \|\nabla u\|_{L^\infty}$ such that

$$\|\nabla u_n\|_{L^\infty} < K', \qquad \forall n \in \mathbb{N}.$$

Consider for instance the Lagrangian $f(\xi) = \exp(-|\xi|)$ presented in Remark 1.10. Then because $f$ is bounded from below by $\exp(-K')$ on $B_{K'}$, necessarily,

$$\liminf_{n \to +\infty} \int_\Omega f(\nabla u_n(x))dx \geq |\Omega|e^{-K'} > 0 = E[f^{**}](u).$$

In particular, (1.17) cannot hold. More generally, if $f : \Omega \times \mathbb{R} \times \mathbb{R}^N \to [0, \infty)$ is any Lagrangian satisfying Hypothesis $\mathcal{A}$, by (1.14),

$$\int_\Omega (\tilde{f}_{K'})^{**}(x, u(x), \nabla u(x))dx \leq \liminf_{n \to +\infty} \int_\Omega f(x, u_n(x), \nabla u_n(x))dx.$$

Thus, recalling that $(\tilde{f}_{K'})^{**} \geq f^{**}$, the existence of a $K' \geq \|\nabla u\|_{L^\infty}$ such that

$$(\tilde{f}_{K'})^{**}(x, u(x), \nabla u(x)) = f^{**}(x, u(x), \nabla u(x)), \qquad \text{for a.e. } x \in \Omega,$$

is a necessary condition in order to reach the infimum in (1.15). Actually this is also a sufficient condition, as it is stated in the following Theorem 1.22.

The previous discussion allows us to characterize the property (named condition ($\mathcal{K}$) in the sequel) which is needed on the Lagrangian to ensure the existence of a sequence satisfying (1.16) and (1.17).

**Definition 1.20.** We say that $f : \Omega \times \mathbb{R} \times \mathbb{R}^N \to [0, \infty]$ satisfies condition ($\mathcal{K}$) if for every compact interval $I \subset \mathbb{R}$ and $K > 0$, there exists $K' \geq K$ and $\tilde{\Omega} \subset \Omega$ such that $|\Omega \setminus \tilde{\Omega}| = 0$ and

$$(\tilde{f}_{K'})^{**} = f^{**} \quad \text{on} \quad \tilde{\Omega} \times I \times B_K. \tag{$\mathcal{K}$}$$

**Remark 1.21.** A simple observation is that, if $g : \mathbb{R}^N \to [0, \infty]$, then

$$(\tilde{g}_{K'})^{**} \geq (\tilde{g}_{K''})^{**} \geq g^{**} \quad \text{on} \quad B_K,$$

for every $K'' \geq K' \geq K > 0$ (by Remark 1.10). Therefore, if condition ($\mathcal{K}$) applies for $f$ and some $K' \geq K$, then it will still apply for any $K'' \geq K'$.

**Theorem 1.22.** *Let $f : \Omega \times \mathbb{R} \times \mathbb{R}^N \to [0, \infty)$ satisfy Hypothesis $\mathcal{A}$. Assume that $f$ also satisfies condition ($\mathcal{K}$). Then for every $u \in W^{1,\infty}(\Omega)$, there exists a sequence $(u_n) \subset W^{1,\infty}_u(\Omega)$ such that*

$$u_n \rightharpoonup^* u \quad \text{in} \quad W^{1,\infty}(\Omega),$$

*and*

$$\lim_{n \to +\infty} \int_\Omega f(x, u_n(x), \nabla u_n(x))dx = \int_\Omega f^{**}(x, u(x), \nabla u(x))dx.$$

*In particular, the sequence $(u_n)$ reaches the minimum for (1.15).*

*Proof.* Let $u$ be in $W^{1,\infty}(\Omega)$, then there exists a compact interval $I \subset \mathbb{R}$ such that
$$u(x) \in I,$$
for almost every $x$ in $\Omega$ and let $K > 0$ be such that
$$\|\nabla u\|_{L^\infty} < K.$$
By assumption we can take $K' \geq K$ such that condition $(\mathcal{K})$ applies, and so in particular
$$(\tilde{f}_{K'})^{**}(x, u(x), \nabla u(x)) = f^{**}(x, u(x), \nabla u(x)), \quad \text{for a.e. } x \in \Omega.$$
By Proposition 1.17 there exists a sequence $(u_n) \subset W_u^{1,\infty}(\Omega)$ such that
$$u_n \rightharpoonup^* u \quad \text{in} \quad W^{1,\infty}(\Omega),$$
and
$$\int_\Omega f(x, u_n(x), \nabla u_n(x))dx \to \int_\Omega (\tilde{f}_{K'})^{**}(x, u(x), \nabla u(x))dx = \int_\Omega f^{**}(x, u(x), \nabla u(x))dx.$$
Since by Theorem 1.19,
$$\inf \left\{ \liminf_{n \to +\infty} \int_\Omega f(x, u_n(x), \nabla u_n(x))dx \;\middle|\; \begin{array}{c} (u_n) \subset W_u^{1,\infty}(\Omega) \\ u_n \rightharpoonup^* u \end{array} \right\} = \int_\Omega f^{**}(x, u(x), \nabla u(x))dx.$$
then $(u_n)$ reaches the minimum in (1.15). $\square$

## 1.4 Sufficient geometric and analytic conditions for $(\mathcal{K})$.

Now we seek conditions on $f$ so that the Lagrangian satisfies condition $(\mathcal{K})$. We begin by stating some results when $f$ depends only on the third variable $\xi$.

**Lemma 1.23.** *Let $f : \mathbb{R}^N \to [0, \infty]$. Assume that $A \subset B \subset \mathbb{R}^N$ satisfy that for any $\xi_A \in A$ and $\xi_B \notin B$, there exists $\xi \in [\xi_A, \xi_B] \cap B$ such that $f(\xi) = f^{**}(\xi)$. Then*
$$(\tilde{f}_B)^{**} = f^{**} \quad on \quad A.$$

*Proof.* Firstly, one has $\tilde{f}_B \geq f$, hence $(\tilde{f}_B)^{**} \geq f^{**}$. We turn to the proof of the converse inequality: assume by contradiction that it is false, then without loss of generality, one may assume that $0 \in A$ and that $(\tilde{f}_B)^{**}(0) > f^{**}(0)$. In particular, thanks to Proposition 1.9, there must exist an affine map $\ell$ on $\mathbb{R}^N$ such that
$$\ell \leq \tilde{f}_B \quad \text{and} \quad \ell(0) > f^{**}(0). \tag{1.18}$$
Define $g := f^{**} - \ell$. Then $g$ is a convex function satisfying $g(0) < 0$ and
$$g(\zeta) \geq 0, \quad \forall \zeta \in \mathbb{R}^N \setminus B. \tag{1.19}$$
Indeed, for $\zeta \notin B$, there exists by assumption $\xi \in [0, \zeta] \cap B$ such that $f(\xi) = f^{**}(\xi)$. If $f(\xi) = +\infty$, then $g(\xi) = +\infty$. Otherwise,
$$g(\xi) = f^{**}(\xi) - \ell(\xi) \geq f^{**}(\xi) - \tilde{f}_B(\xi) = f^{**}(\xi) - f(\xi) = 0.$$
Therefore, in any case $g(\xi) \geq 0$ and because $g(0) < 0$, (1.19) holds. We derive from (1.19) that $\ell \leq f$ on $\mathbb{R}^N$, indeed:

- on $B$, one has $f = \tilde{f}_B$ and the conclusion follows from (1.18) in this case;
- on $\mathbb{R}^N \setminus B$, this comes from (1.19) and the fact that $f^{**} \leq f$.

Finally, one gets that $\ell$ is an affine map which bounds $f$ from below, thus $\ell \leq f^{**}$. But this is in clear contradiction with (1.18), therefore the proof is complete. □

Lemma 1.23 has a few consequences:

**Corollary 1.24.** *Let $f : \mathbb{R}^N \to [0, \infty]$.*

- *If $A \subset \mathbb{R}^N$ is closed and $f = f^{**}$ on $\partial A$, then $(\tilde{f}_A)^{**} = f^{**}$ on $A$.*

- *Assume that the detachment set of $f$ defined by*

$$\mathcal{D}(f) := \{\xi \in \mathbb{R}^N \ : \ f^{**}(\xi) < f(\xi)\}$$

*has uniformly bounded connected components (i.e. there exists $M > 0$ such that $\operatorname{diam}(C) \leq M$ for any such component). Then for any $K > 0$, there exists $K' \geq K$ such that*

$$(\tilde{f}_{K'})^{**} = f^{**} \quad on \quad B_K.$$

*Proof.*
- Taking $B = A$ in Lemma 1.23 and noticing that by assumption $\partial A \subset A$, it gives the result.

- Take $K > 0$ and let $K' > K + M$. Then for any $\xi_1 \in B_K$, $\xi_2 \notin B_{K'}$, one has $|\xi_1 - \xi_2| > M$ and as such, there must exist some $\xi \in [\xi_1, \xi_2]$ with $|\xi| \leq K'$ such that $f(\xi) = f^{**}(\xi)$. The result follows from Lemma 1.23.

□

**Remark 1.25.** The assumption that $A$ is closed in the first assertion of the previous corollary is crucial. Take for instance $f : \mathbb{R} \to [0, \infty)$ defined by

$$f(\xi) = \begin{cases} |\xi|^2 & \text{if} \quad \xi < 1 \\ 0 & \text{otherwise.} \end{cases}$$

One can verify easily (using the second point of Proposition 1.9 for instance) that

$$f^{**}(\xi) = \begin{cases} |\xi|^2 & \text{if} \quad \xi \leq 0 \\ 0 & \text{otherwise.} \end{cases}$$

Taking $A = (-\infty, 1)$, then $f = f^{**}$ on $\partial A = \{1\}$. However, $f$ is convex on $A$ and therefore

$$(\tilde{f}_A)^{**} = f \neq f^{**} \quad on \quad A.$$

We now generalize Corollary 1.24 to the case of $x$ and $u$ dependency, whose proof is immediate:

**Corollary 1.26.** *Let $f : \Omega \times \mathbb{R} \times \mathbb{R}^N \to [0, \infty]$. Let $I$ be a compact interval of $\mathbb{R}$.*

- *If there exists $K > 0$ such that*

$$f(x, u, \xi) = f^{**}(x, u, \xi), \qquad \text{for a.e. } x \in \Omega, \ \forall u \in I, \ \forall \xi \in \partial B_K, \qquad (1.20)$$

*Then*

$$(\tilde{f}_K)^{**}(x, u, \xi) = f^{**}(x, u, \xi), \qquad \text{for a.e. } x \in \Omega, \ \forall u \in I, \ \forall \xi \in B_K.$$

*In particular, if (1.20) holds for a sequence $(K_n)$ of positive numbers such that $K_n \to +\infty$, and for every such $I$, then $f$ satisfies $(\mathcal{K})$.*

- *Assume that the connected components of the sets*

$$\mathcal{D}(f, x, u) := \{\xi \in \mathbb{R}^N \ : \ f^{**}(x, u, \xi) < f(x, u, \xi)\}, \qquad x \in \Omega, \ u \in I,$$

*are uniformly bounded. That is, there exists $M > 0$ such that for a.e. $x \in \Omega$, for all $u \in I$, every connected component of $\mathcal{D}(f, x, u)$ has diameter smaller than $M$. If this holds for all compact interval $I$, then $f$ satisfies $(\mathcal{K})$.*

The second condition that we present here, independent of the geometry of the detachment set, is the superlinearity of the Lagrangian with respect to the variable $\xi$. In the following result, we will consider the dependency in $(x, u)$ as auxiliary, taken into account in the parametric set $\Gamma$. It is an abstract statement and one may identify $\Gamma$ as $\Omega \times I$, for $I \subset \mathbb{R}$ a compact interval (see Corollary 1.30 below which gives a more explicit statement).

**Theorem 1.27.** *Let $\Gamma$ be a set and $f : \Gamma \times \mathbb{R}^N \to [0, \infty]$ be such that:*

- *there exists $\Phi : \mathbb{R}^N \to [0, \infty)$ superlinear such that*

$$f(s, \xi) \geq \Phi(\xi), \qquad \forall s \in \Gamma, \ \forall \xi \in \mathbb{R}^N; \qquad (1.21)$$

- *for any $\rho > 0$, there exists $\rho' \geq \rho$ such that $(\tilde{f}_{\rho'})^{**}$ is bounded on $\Gamma \times B_\rho$.*

*Then for any $K > 0$, there exists $K' \geq K$ such that*

$$(\tilde{f}_{K'})^{**} = f^{**} \quad \text{on} \quad \Gamma \times B_K.$$

**Remark 1.28.** Following the same convention introduced in Definition 1.12, the convexification in this Theorem is to be understood with respect to the variable $\xi \in \mathbb{R}^N$.

*Proof.* Let $\rho' > 0$ satisfy the second assumption for $\rho = K + 1$, and define

$$M = \sup_{\Gamma \times B_{K+1}} (\tilde{f}_{\rho'})^{**} < +\infty.$$

By (1.21), there exists $K' \geq \rho'$ such that

$$f(s, \xi) \geq (M + 1)|\xi| + M, \qquad \forall s \in \Gamma, \ \forall \xi \in \mathbb{R}^N \setminus B_{K'}. \qquad (1.22)$$

Fix now $\varepsilon \in (0, 1]$ and consider $(s_0, \xi_0) \in \Gamma \times B_K$, then by Proposition 1.9, there exists an affine map $\ell$ of slope $\zeta \in \mathbb{R}^N$, (i.e. $\nabla \ell \equiv \zeta$), such that

$$\ell \leq \tilde{f}_{K'}(s_0, \cdot) \quad \text{and} \quad \ell(\xi_0) \geq (\tilde{f}_{K'})^{**}(s_0, \xi_0) - \varepsilon. \tag{1.23}$$

Owing to (1.23), one has

$$-1 \leq (\tilde{f}_{K'})^{**}(s_0, \xi_0) - \varepsilon \leq \ell(\xi_0) \leq \ell\left(K\frac{\zeta}{|\zeta|}\right) = \ell\left((K+1)\frac{\zeta}{|\zeta|}\right) - |\zeta| \leq M - |\zeta|,$$

because $\ell \leq (\tilde{f}_{K'})^{**}(s_0, \cdot)$ and $\sup_{B_{K+1}} (\tilde{f}_{K'})^{**}(s_0, \cdot) \leq M$. Hence $|\zeta| \leq M + 1$. Our goal is to show that $\ell$ is an affine minorant of $f(s_0, \cdot)$. The fact that $\ell \leq f(s_0, \cdot)$ on $B_{K'}$ is an immediate consequence of (1.23). Now let us consider $\xi$ which is not in $B_{K'}$, then by (1.22):

$$\ell(\xi) = \langle \xi, \zeta \rangle + \ell(0) \leq (M+1)|\xi| + M \leq f(s_0, \xi),$$

which proves that $\ell \leq f(s_0, \cdot)$ on $\mathbb{R}^N$ and thus that $\ell \leq f^{**}(s_0, \cdot)$. In particular, one has by (1.23) :

$$(\tilde{f}_{K'})^{**}(s_0, \xi_0) \leq f^{**}(s_0, \xi_0) + \varepsilon,$$

letting $\varepsilon \to 0$, it gives $(f_{K'})^{**} \leq f^{**}$ on $\Gamma \times B_K$. The converse inequality being obvious, we finally obtain that

$$(f_{K'})^{**} = f^{**} \quad \text{on} \quad \Gamma \times B_K.$$

$\square$

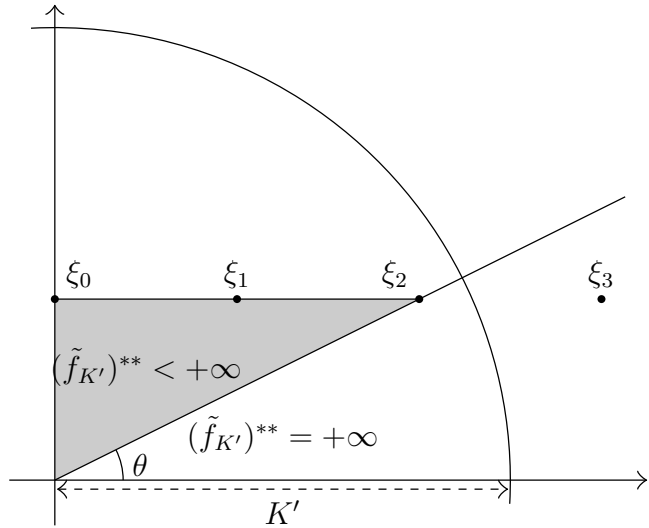

Figure 1: $(\tilde{f}_{K'})^{**}$ is infinite on the angular sector guided by $\theta$. Yet, it is not the case for $f^{**}$.

**Remark 1.29.** The second assumption of Theorem 1.27 is crucial in this result. So much that it may fail without it: consider for instance the function (depending only on $\xi$), $f$ :

$\mathbb{R}^2 \to [0, \infty]$ defined by:

$$\begin{cases} f(0,0) &= 0 \\ f(\xi_n) &= |\xi_n|^2 & \forall n \in \mathbb{N} \\ f(\xi) &= +\infty & \text{otherwise,} \end{cases}$$

where $\xi_n := (n, 1)$. Then $f$ is superlinear (taking for instance $\Phi(\xi) = |\xi|^2$ is enough), but we claim that it does not satisfy $(\mathcal{K})$. Indeed, $f$ is finite on the set

$$A := \{(0,0)\} \cup \{\xi_0, \xi_1, \dots\},$$

therefore $f^{**}$ is finite on $\operatorname{conv} A \supset [0, \infty[\times]0, 1]$. But for any $K' > 0$, $B_{K'}$ contains only a finite number of the $\xi_n$ and thus $(\tilde{f}_{K'})^{**}$ is infinite on an angular sector (see Figure 1).

For the sake of comprehensibility, we now state a version of Theorem 1.27 in the case where $f : \Omega \times \mathbb{R} \times \mathbb{R}^N \to [0, \infty]$:

**Corollary 1.30.** *Let $f : \Omega \times \mathbb{R} \times \mathbb{R}^N \to [0, \infty]$. Assume that for every compact interval $I \subset \mathbb{R}$,*

- *there exists $\Phi_I : \mathbb{R}^N \to [0, \infty)$ superlinear such that*

$$f(x, u, \xi) \geq \Phi_I(\xi), \qquad \text{for a.e. } x \in \Omega, \ \forall u \in I, \ \forall \xi \in \mathbb{R}^N; \qquad (1.24)$$

- *there exists $\tilde{\Omega} \subset \Omega$ such that $|\Omega \setminus \tilde{\Omega}| = 0$ and for any $\rho > 0$, $f$ is bounded on $\tilde{\Omega} \times I \times B_\rho$.*

*Then $f$ satisfies condition $(\mathcal{K})$.*

*Proof.* Let $\Omega_0 \subset \Omega$ a measurable set with full measure such that (1.24) holds for every $x \in \Omega_0$. Let $\Gamma = (\Omega_0 \cap \tilde{\Omega}) \times I$. Then the first condition of Theorem 1.27 holds by construction, and the second comes from the fact that

$$(\tilde{f}_\rho)^{**} \leq f, \quad \text{on} \quad \Omega \times I \times B_\rho,$$

thus the choice $\rho' = \rho$ is enough. $\square$

On another note, Theorem 1.27 allows us to prove Theorem 1.16, which gave a sufficient condition for $f^{**}$ to be Carathéodory..

*Proof of Theorem 1.16.* By Lemma 1.13, we know that $f^{**}$ is Lebesgue-Borel measurable. It is therefore enough to show that it is continuous with respect to $(u, \xi)$ for a.e. $x \in \Omega$. Let $A_0$ be a measurable subset of $\Omega$ with full measure such that (1.10) holds for every $x \in A_0$ and every compact interval $I$. According to Hypothesis $\mathcal{A}$-c), for any integer $n \geq 1$, there exists $a_n \in L^1(\Omega, [0, \infty))$ and a measurable set $A_n \subset \Omega$ with full measure such that

$$f(x, u, \xi) \leq a_n(x) < +\infty, \qquad \forall x \in A_n, \ \forall (u, \xi) \in [-n, n] \times B_n. \qquad (1.25)$$

For $x \in A := \cap_{n\geq 0} A_n$, define $g : (u, \xi) \mapsto f(x, u, \xi)$. We will show that $g$ satisfies the assumption of Theorem 1.27 with $\Gamma = I$. Fix $I \subset \mathbb{R}$ a compact interval. Firstly by (1.25), $g$ is bounded on bounded sets and thus the second assumption of Theorem 1.27 is satisfied with $\rho' = \rho$. For the first assumption, notice that by (1.10), it is satisfied with $\Phi := \Phi_{x,I}$. We can therefore finally apply Theorem 1.27, which tells us that $g$ satisfies (1.9) of Proposition 1.15 and thus that $g^{**}$ is continuous. The conclusion follows from the fact that $A$ has full measure in $\Omega$, and hence we have proved that $(u, \xi) \mapsto f^{**}(x, u, \xi)$ is continuous for a.e. $x \in \Omega$: $f^{**}$ is a Carathéodory function. $\square$

**Remark 1.31.** In this Section 1.4, we have given two distinct sets of assumptions ensuring that a Lagrangian satisfies $(\mathcal{K})$. On the one hand, a boundedness condition on the detachment set (Corollaries 1.24 and 1.26). On the other hand, the superlinearity of the Lagrangian (Theorem 1.27 and Corollary 1.30). We wish to show that these sets of assumptions are independent and none implies the other.

- Let $f : (\xi_1, \xi_2) \mapsto (1 + \sin(\xi_1))(1 + \sin(\xi_2))$. Then $f$ is clearly not superlinear, however $f^{**} \equiv 0$ and thus $f$ satisfies the detachment set condition, (second point of Corollary 1.24) because $f = f^{**}$ on a square grid of side-length $2\pi$.

- Conversely, let $f : \mathbb{R}^2 \to [0, \infty)$ defined by

$$f(\xi_1, \xi_2) := \big(|\xi_1| - 1\big)^2 + |\xi_2|^2.$$

It is clearly superlinear. However, we show that the line $\{\xi_1 = 0\}$ is contained in the detachment set of $f$, thus contradicting the boundedness condition. Indeed, notice that by Proposition 1.9, for every $\xi_2 \in \mathbb{R}$,

$$\begin{aligned}
f^{**}(0, \xi_2) &\leq \frac{1}{2}(f(-1, \xi_2) + f(1, \xi_2)) \\
&= |\xi_2|^2 \\
&< |\xi_2|^2 + 1 \\
&= f(0, \xi_2).
\end{aligned}$$

# 2  Approximation in $W^{1,p}(\Omega)$ for the autonomous case

In this paragraph, we state an important result which answers the question of the Lavrentiev gap when the Lagrangian is *autonomous*, that is $f = f(u, \xi)$ does not depend on $x$. Some regularity is still needed (in particular, the convexity with respect to $\xi$). The detailed proof may be found in [8].

**Theorem 2.1.** *Let $f : \mathbb{R} \times \mathbb{R}^N \to [0, \infty)$ be a continuous Lagrangian, which is convex in $\xi$ (the last variable), and let $\varphi \in W^{1,\infty}(\Omega)$. Then for any $u \in W^{1,1}_\varphi(\Omega)$, there exists a sequence $(u_n) \subset W^{1,\infty}_\varphi(\Omega)$ such that*

$$u_n \to u \quad \text{strongly in} \quad W^{1,1}(\Omega),$$

*and*

$$\lim_{n\to+\infty} \int_\Omega f(u_n(x), \nabla u_n(x))dx = \int_\Omega f(u(x), \nabla u(x))dx.$$

*Moreover, if $u \in W^{1,p}(\Omega)$ for some $p > 1$, then the sequence $(u_n)$ can be chosen so that $u_n \to u$ strongly in $W^{1,p}(\Omega)$.*

**Remark 2.2.** We can see Theorem 2.1 through another perspective by defining another functional:

$$E_{rel}[f](u) = \inf \left\{ \liminf_{n\to+\infty} \int_\Omega f(u_n(x), \nabla u_n(x))dx \;\middle|\; \begin{array}{c} (u_n) \subset W^{1,\infty}_\varphi(\Omega) \\ u_n \to_{W^{1,1}} u \end{array} \right\}.$$

By the Fatou Lemma we have that if $f(u, \xi)$ is continuous then for every $u \in W^{1,1}_\varphi(\Omega)$

$$E[f](u) \leq E_{rel}[f](u),$$

with an equality for every $u \in W^{1,\infty}_\varphi(\Omega)$. Now Theorem 2.1 implies that if $g(u, \xi)$ is continuous and convex with respect to the last variable then

$$E[g] = E_{rel}[g] \quad \text{on} \quad W^{1,1}_\varphi(\Omega).$$

Now a natural question arising from this statement is to wonder whether the assumption concerning convexity is needed for this result to hold. In this Section, we will extend this result to a larger class of Lagrangians. In the autonomous case we have a simpler set of assumptions.

**Hypothesis $\mathcal{B}$.** The function $f : \mathbb{R} \times \mathbb{R}^N \to [0, \infty)$ satisfies

a) $f$ is Borel,

b) $f$ is continuous with respect to $u$ uniformly as $\xi$ varies in bounded sets of $\mathbb{R}^N$,

c) $f$ is bounded on bounded sets,

d) $f^{**}$, the bipolar of $f$ with respect to the second variable, is continuous in $(u, \xi)$.

In the autonomous case Hypothesis $\mathcal{B}$ implies Hypothesis $\mathcal{A}$ (this is pretty clear but a detailed proof is given in [6, Lemma 18]). We recall that Theorem 1.16 gives a sufficient condition for $\mathcal{B}$-d) to hold: it is enough to assume that $f$ is superlinear with respect to $\xi$ uniformly as $u$ varies in bounded sets of $\mathbb{R}$. It might be of interest to note that again, Hypothesis $\mathcal{B}$ does not request that $f$ is continuous with respect to $(u, \xi)$.

Of course, Theorem 2.1 implies the non occurrence of the Lavrentiev phenomenon for continuous and autonomous Lagrangians which are convex in $\xi$. Actually, more is true, as stated in the following result:

**Theorem 2.3.** *Let $f : \mathbb{R} \times \mathbb{R}^N \to [0, \infty)$ satisfy Hypothesis $\mathcal{B}$. Then no Lavrentiev phenomenon occurs for the integral functional $E[f]$, that is:*

$$\inf \left\{ E[f](u) \;:\; u \in W^{1,1}_\varphi(\Omega) \right\} = \inf \left\{ E[f](u) \;:\; u \in W^{1,\infty}_\varphi(\Omega) \right\}.$$

*Proof.* For the sake of simplicity, write $X = W^{1,1}_\varphi(\Omega)$ and $Y = W^{1,\infty}_\varphi(\Omega)$. Then

$$\inf_X E[f] \geq \inf_X E[f^{**}] = \inf_Y E[f^{**}] = \inf_Y E[f] \geq \inf_X E[f].$$

Indeed, the first and last inequalities are clear, the second comes from Theorem 2.1 and the third from Theorem 1.19. □

As a byproduct of this proof, we deduce that (if $f$ satisfies Hypothesis $\mathcal{B}$), for any $p \in [1, \infty)$,

$$\inf \left\{ E[f](u) \ : \ u \in W^{1,p}_\varphi(\Omega) \right\} = \inf \left\{ E[f^{**}](u) \ : \ u \in W^{1,p}_\varphi(\Omega) \right\}. \tag{2.1}$$

## 2.1 Weak relaxation on $W^{1,1}(\Omega)$.

Now we extend Theorem 2.1 to the non convex case. The following theorem is a generalization of [6, Theorem 30].

**Theorem 2.4.** *Let $f : \mathbb{R} \times \mathbb{R}^N \to [0, \infty)$ satisfy Hypothesis $\mathcal{B}$, and let $\varphi \in W^{1,\infty}(\Omega)$. Then for any $u \in W^{1,1}_\varphi(\Omega)$, there exists a sequence $(u_n) \subset W^{1,\infty}_\varphi(\Omega)$ such that*

$$u_n \to u \quad in \quad L^1(\Omega),$$

*and*

$$\lim_{n \to +\infty} \int_\Omega f(u_n(x), \nabla u_n(x)) dx = \int_\Omega f^{**}(u(x), \nabla u(x)) dx.$$

*Moreover,*

- *if there exists $\Phi : \mathbb{R}^N \to [0, \infty)$ superlinear such that*

$$f(u, \xi) \geq \Phi(\xi), \qquad \forall u \in \mathbb{R}, \ \forall \xi \in \mathbb{R}^N, \tag{2.2}$$

  *then the sequence $(u_n)$ can be chosen so that $u_n \rightharpoonup u$ weakly in $W^{1,1}(\Omega)$;*

- *if for some $p \in (1, \infty)$ it holds that*

$$f(u, \xi) \geq c_1 |\xi|^p - c_2, \qquad \forall u \in \mathbb{R}, \ \forall \xi \in \mathbb{R}^N, \tag{2.3}$$

  *for some $c_1, c_2 > 0$, and $u \in W^{1,p}(\Omega)$ then the sequence $(u_n)$ can be taken so that $u_n \rightharpoonup u$ weakly in $W^{1,p}(\Omega)$.*

**Remark 2.5.** Notice that, compared to Theorem 2.1, the assumptions on $f$ are relaxed: indeed in particular, $f$ is not assumed to be convex in $\xi$. However, the conclusion is slightly weaker in the sense that the strong convergence of Theorem 2.1 was replaced by an $L^1$ (or $W^{1,1}$-weak) convergence. Apart from that, both results have similar conclusions, up to noticing that $f = f^{**}$ in the assumptions of Theorem 2.1.

*Proof of Theorem 2.4.* Fix $u \in W^{1,1}_\varphi(\Omega)$. If $\int_\Omega f^{**}(u(x), \nabla u(x))dx = +\infty$, the result follows from Fatou Lemma. Indeed, take any sequence $(u_n) \subset W^{1,\infty}_\varphi(\Omega)$ converging to $u$ strongly in $W^{1,1}(\Omega)$ (such a sequence always exists). Then up to a subsequence, $(u_n(x), \nabla u_n(x)) \to (u(x), \nabla u(x))$ for a.e. $x \in \Omega$, and by continuity of $f^{**}$, it holds that $f^{**}(u_n(x), \nabla u_n(x)) \to f^{**}(u(x), \nabla u(x))$ for a.e. $x \in \Omega$. Thus

$$\liminf_{n \to +\infty} \int_\Omega f(u_n(x), \nabla u_n(x))dx \geq \liminf_{n \to +\infty} \int_\Omega f^{**}(u_n(x), \nabla u_n(x))dx$$
$$\geq \int_\Omega f^{**}(u(x), \nabla u(x))dx$$
$$= +\infty.$$

We assume from now on that $E[f^{**}](u) < +\infty$. By Theorem 2.1, there exists a sequence $(v_n) \subset W^{1,\infty}_\varphi(\Omega)$ such that $v_n \to u$ in $W^{1,1}(\Omega)$ and $E[f^{**}](v_n) \to E[f^{**}](u)$. Now for each $n$, one may use Theorem 1.19 to find a sequence $(v_n^k) \subset W^{1,\infty}_\varphi(\Omega)$ such that $\|v_n^k - v_n\|_{L^\infty} \to 0$ (in particular, $v_n^k \to v_n$ in $L^1$) as $k \to +\infty$ and

$$E[f](v_n^k) \xrightarrow{k \to +\infty} E[f^{**}](v_n).$$

Using a diagonal argument, one obtains a sequence $(u_n) \subset W^{1,\infty}_\varphi(\Omega)$ such that $u_n \to u$ in $L^1$ and $E[f](u_n) \to E[f^{**}](u)$. To conclude, notice that conditions (2.2) and (2.3) respectively gives weak compactness for the approximating sequence in $W^{1,1}$ (by the Dunford-Pettis Theorem) and $W^{1,p}$. Therefore it is an immediate corollary that the approximating sequence weakly converges to $u$ in that case. □

**Remark 2.6.** We can see Theorem 2.4 from a functional point of view.
We define

$$\tilde{E}_1[f](u) := \inf \left\{ \liminf_{n \to +\infty} \int_\Omega f(u_n(x), \nabla u_n(x))dx \ \middle| \ \begin{array}{c} (u_n) \subset W^{1,\infty}_\varphi(\Omega) \\ u_n \rightharpoonup_{W^{1,1}} u \end{array} \right\}. \quad (2.4)$$

If $f^{**}$ is continuous in $(u, \xi)$, by [36, Chapter 4, Theorem 4.5] the integral functional

$$u \mapsto \int_\Omega f^{**}(u(x), \nabla u(x))dx,$$

is weakly lower semicontinuous on $W^{1,1}(\Omega)$. So we have

$$\int_\Omega f^{**}(u(x), \nabla u(x))dx \leq \liminf_{n \to +\infty} \int_\Omega f(v_n(x), \nabla v_n(x))dx,$$

for every $(v_n) \subset W^{1,1}(\Omega)$ such that $v_n \rightharpoonup_{W^{1,1}} u$, that is

$$\int_\Omega f^{**}(u(x), \nabla u(x))dx \leq \tilde{E}_1[f](u).$$

Now we know by Theorem 2.4 that if $f$ satisfies Hypothesis $\mathcal{B}$ and is uniformly superlinear then for every $u \in W^{1,1}_\varphi(\Omega)$,

$$E[f^{**}](u) = \tilde{E}_1[f](u),$$

and furthermore the infimum in (2.4) is actually a minimum.

## 2.2 A geometric condition to recover strong convergence.

In Theorem 2.4 we do not have the strong convergence of the approximating sequence. Now we study other cases in which the approximating sequence converges strongly in $W^{1,1}(\Omega)$. Before stating the main result, we will need a useful lemma.

**Lemma 2.7.** *Let $v \in L^1(\Omega, \mathbb{R}^N)$. Then there exists $\Phi : \mathbb{R}^N \to [0, \infty)$ convex and superlinear such that*
$$\int_\Omega \Phi(v(x))dx < +\infty.$$

*Proof.* The proof is rather straightforward and uses a similar argument to the one used in the proof of the Dunford-Pettis Theorem. Let $(M_n)$ be an increasing sequence of positive real numbers such that $M_n \to +\infty$ and
$$\int_{\{|v| \geq M_n\}} |v|dx \leq 2^{-n}, \qquad \forall n \in \mathbb{N}. \tag{2.5}$$

Then define
$$\Phi(\xi) := \sum_{n \in \mathbb{N}} \left(|\xi| - M_n\right)^+, \qquad \forall \xi \in \mathbb{R}^N,$$

where we used the standard notation $a^+ := \max(a, 0)$. The assumption that $M_n \to +\infty$ ensures that the sum in the definition of $\Phi$ is effectively finite at each point, and therefore $\Phi$ takes values into $[0, \infty)$. Each of the functions $\xi \mapsto (|\xi| - M_n)^+$ is convex, thus so is $\Phi$ as a supremum of convex functions. Also, the identity $\Phi(\xi)/|\xi| = \sum_n (1 - M_n/|\xi|)^+$ ensures that $\Phi$ is superlinear. To achieve the proof, it is finally enough to use (2.5) and notice the following:
$$\int_\Omega \Phi(v(x))dx = \sum_{n=0}^{+\infty} \int_{\{|v| \geq M_n\}} \left(|v| - M_n\right)dx \leq \sum_{n=0}^{+\infty} 2^{-n} < +\infty.$$
$\square$

**Theorem 2.8.** *Let $f : \mathbb{R} \times \mathbb{R}^N \to [0, \infty)$ satisfies Hypothesis B, and let $\varphi \in W^{1,\infty}(\Omega)$. Assume $u \in W^{1,p}_\varphi(\Omega)$ for some $p \in [1, \infty)$, and satisfies*
$$\int_\Omega f^{**}(u(x), \nabla u(x))dx = \int_\Omega f(u(x), \nabla u(x))dx. \tag{2.6}$$

*Then there exists a sequence $(u_n) \subset W^{1,\infty}_\varphi(\Omega)$ such that*
$$u_n \to u \quad \text{strongly in} \quad W^{1,p}(\Omega),$$

*and*
$$\lim_{n \to +\infty} \int_\Omega f(u_n(x), \nabla u_n(x))dx = \int_\Omega f(u(x), \nabla u(x))dx.$$

We stated this result as being geometric in the following sense: 2.6 means that the map $(u, \nabla u)$ does not see the detachment set of $f$ on $\mathbb{R} \times \mathbb{R}^N$ (see Remark 2.10 below). In fact, equality (2.6) happens to also be a necessary for the existence of such an approximating sequence (see Remark 2.12).

*Proof.* Assume first that $p = 1$. In the case where $E[f](u) = +\infty$, this follows from the Fatou lemma just as in the proof of Theorem 2.4. Otherwise, we assume $E[f](u) < +\infty$. By Lemma 2.7, there exists $\Phi : \mathbb{R}^N \to [0, \infty)$ convex and superlinear such that

$$\int_\Omega \Phi(\nabla u(x))dx < +\infty.$$

We define

$$g(u, \xi) := f(u, \xi) + \Phi(\xi) + \sqrt{1 + |\xi|^2}.$$

Because $g(u, \xi) \geq \Phi(\xi)$, Theorem 1.16 ensures that $g^{**}$ is continuous and thus $g$ satisfies Hypothesis $\mathcal{B}$. By Theorem 2.4, there exists a sequence $(u_n) \subset W^{1,\infty}_\varphi(\Omega)$ such that $u_n \rightharpoonup u$ weakly in $W^{1,1}(\Omega)$ and

$$\int_\Omega g(u_n(x), \nabla u_n(x))dx \to \int_\Omega g^{**}(u(x), \nabla u(x))dx$$
$$\leq E[f](u) + \int_\Omega \Phi(\nabla u(x))dx + \int_\Omega \sqrt{1 + |\nabla u(x)|^2}dx.$$

(The inequality is a simple consequence of $g^{**} \leq g$.) Therefore,

$$\limsup_{n \to +\infty} \left( E[f](u_n) + \int_\Omega \Phi(\nabla u_n(x))dx + \int_\Omega \sqrt{1 + |\nabla u_n(x)|^2}dx \right) \quad (2.7)$$
$$\leq E[f](u) + \int_\Omega \Phi(\nabla u(x))dx + \int_\Omega \sqrt{1 + |\nabla u(x)|^2}dx < +\infty$$

By convexity with respect to the last variable and continuity, we have again by Tonelli's result (see Lemma 1.18) the following estimates

$$E[f](u) = \int_\Omega f^{**}(u(x), \nabla u(x))dx \leq \liminf_{n \to +\infty} \int_\Omega f^{**}(u_n(x), \nabla u_n(x))dx$$
$$\leq \liminf_{n \to +\infty} E[f](u_n),$$

as well as

$$\int_\Omega \Phi(\nabla u(x))dx \leq \liminf_{n \to +\infty} \int_\Omega \Phi(\nabla u_n(x))dx,$$

$$\int_\Omega \sqrt{1 + |\nabla u(x)|^2}dx \leq \liminf_{n \to +\infty} \int_\Omega \sqrt{1 + |\nabla u_n(x)|^2}dx,$$

and thus, combining this with (2.7), one see that we have in fact convergence of each of these terms. In particular:

$$E[f](u_n) \to E[f](u) \quad \text{and} \quad \int_\Omega \sqrt{1 + |\nabla u_n|^2}dx \to \int_\Omega \sqrt{1 + |\nabla u|^2}dx.$$

From the second point, we deduce that $u_n \to u$ strongly in $W^{1,1}(\Omega)$ (see Lemma 2.9 below), which achieves the proof. When $p > 1$, the proof is similar, as one only needs to take $g(u, \xi) = f(u, \xi) + |\xi|^p$ in this case.

$\square$

**Lemma 2.9.** *Let $p \in [1, \infty)$ and $(v_n) \subset L^p$. Assume that $v_n \rightharpoonup v$ weakly in $L^p$.*

- *If $p = 1$ and $\int \sqrt{1 + |v_n|^2} \to \int \sqrt{1 + |v|^2}$, then $v_n \to v$ strongly in $L^1$.*
- *If $p > 1$ and $\int |v_n|^p \to \int |v|^p$, then $v_n \to v$ strongly in $L^p$.*

*Proof.* For the case $p = 1$, see [35, Section 1.3.4, Proposition 1]. For $p > 1$, this is a consequence of the uniform convexity of the space $L^p$. $\square$

**Remark 2.10.** It might be relevant to emphasize the fact that, using the inequality $f^{**} \leq f$, condition (2.6) in Theorem 2.8 is satisfied if and only if

$$f(u(x), \nabla u(x)) = f^{**}(u(x), \nabla u(x)), \qquad \text{for a.e. } x \in \Omega.$$

This Theorem 2.8 has as a consequence the following interesting corollary, which prevent the Lavrentiev gap for functions $\overline{u}$ which are already known to be minimizers of the energy.

**Corollary 2.11.** *Let $f : \mathbb{R} \times \mathbb{R}^N \to [0, \infty)$ satisfies Hypothesis B, and let $\varphi \in W^{1,\infty}(\Omega)$. If $\overline{u}$ is a minimizer of $E[f]$ on $W^{1,p}_\varphi(\Omega)$, for some $p \in [1, \infty)$, then there exists a sequence $(u_n) \subset W^{1,\infty}_\varphi(\Omega)$ such that*

$$u_n \to \overline{u} \quad \text{strongly in} \quad W^{1,p}(\Omega),$$

*and*

$$\lim_{n \to +\infty} \int_\Omega f(u_n(x), \nabla u_n(x)) dx = \int_\Omega f(\overline{u}(x), \nabla \overline{u}(x)) dx.$$

*Proof.* By Theorem 2.8, it is enough to prove that $\overline{u}$ satisfies (2.6). Thanks to equality (2.1), we obtain

$$\begin{aligned} E[f^{**}](\overline{u}) &\geq \inf \left\{ E[f^{**}](u) \ : \ u \in W^{1,p}_\varphi(\Omega) \right\} \\ &= \inf \left\{ E[f](u) \ : \ u \in W^{1,p}_\varphi(\Omega) \right\} \\ &= E[f](\overline{u}). \end{aligned}$$

The reverse inequality comes from $f^{**} \leq f$. $\square$

**Remark 2.12.** We can see Theorem 2.8 as a generalization of Theorem 2.1 to the case of non convex Lagrangians. Indeed, if $f = f^{**}$, then (2.6) is trivially satisfied for every $u$ and we recover the original statement of Theorem 2.1. In fact more can be said on this matter: Theorem 2.8 states that (2.6) is a sufficient condition for the existence of a sequence $(u_n)$ strongly converging to $u$ satisfying $E[f](u_n) \to E[f^{**}](u)$. We claim that under the assumption that $f$ is continuous, the converse is also true. Indeed, in that case, by the Fatou Lemma,

$$E[f](u) \leq \liminf_{n \to +\infty} E[f](u_n) = E[f^{**}](u),$$

and the fact that $f^{**} \leq f$ gives the reverse inequality.

## 2.3 Convergence results for convex-dominated Lagrangians.

In the following theorem we show another case in which $E[f](u) = E_{rel}[f](u)$ (the notation $E_{rel}[f]$ was introduced in Remark 2.2), with some interesting implications.

**Theorem 2.13.** *Let $f : \mathbb{R} \times \mathbb{R}^N \to [0, \infty)$ continuous, and let $\varphi \in W^{1,\infty}(\Omega)$. Let $u \in W^{1,p}_\varphi(\Omega)$ for some $p \in [1, \infty)$ and assume that there exists $g : \mathbb{R} \times \mathbb{R}^N \to [0, \infty)$ globally continuous and convex with respect to the last variable, such that $f \leq g$ and*

$$\int_\Omega g(u(x), \nabla u(x)) dx < +\infty. \tag{2.8}$$

*Then there exists a sequence $(u_n) \subset W^{1,\infty}_\varphi(\Omega)$ such that*

$$u_n \to u \quad \text{strongly in} \quad W^{1,p}(\Omega),$$

*and*

$$\lim_{n \to +\infty} \int_\Omega f(u_n(x), \nabla u_n(x)) dx = \int_\Omega f(u(x), \nabla u(x)) dx.$$

This result can be seen as a from of dominated convergence theorem. In fact, it is an immediate consequence of the more general following one:

**Theorem 2.14.** *Let $f : \mathbb{R} \times \mathbb{R}^N \to [0, \infty)$ continuous, and let $\varphi \in W^{1,\infty}(\Omega)$. Let $u \in W^{1,p}_\varphi(\Omega)$ for some $p \in [1, \infty)$ and assume that there exists a continuous $g : \mathbb{R} \times \mathbb{R}^N \to [0, \infty)$ satisfying Hypothesis B, such that $f \leq g$, (2.8) is satisfied, as well as*

$$\int_\Omega g^{**}(u(x), \nabla u(x)) dx = \int_\Omega g(u(x), \nabla u(x)) dx.$$

*Then there exists a sequence $(u_n) \subset W^{1,\infty}_\varphi(\Omega)$ such that*

$$u_n \to u \quad \text{strongly in} \quad W^{1,p}(\Omega),$$

*and*

$$\lim_{n \to +\infty} \int_\Omega f(u_n(x), \nabla u_n(x)) dx = \int_\Omega f(u(x), \nabla u(x)) dx.$$

*Proof.* We begin by noticing that $g$ satisfies the assumptions of Theorem 2.8. Therefore, there exists $(u_n) \subset W^{1,\infty}_\varphi(\Omega)$ such that

$$u_n \to u \quad \text{strongly in} \quad W^{1,p}(\Omega),$$

and

$$\lim_{n \to +\infty} \int_\Omega g(u_n(x), \nabla u_n(x)) dx = \int_\Omega g(u(x), \nabla u(x)) dx. \tag{2.9}$$

Since $f$ and $g - f$ are continuous and non-negative, by the Fatou Lemma,

$$\liminf_{n \to +\infty} E[f](u_n) \geq E[f](u), \tag{2.10}$$

and
$$\liminf_{n \to +\infty} E[g - f](u_n) \geq E[g - f](u). \tag{2.11}$$

Thus, using (2.9) and (2.11),

$$\begin{aligned}
\limsup_{n \to +\infty} E[f](u_n) &= \limsup_{n \to +\infty} \left[ E[g](u_n) - E[g - f](u_n) \right] \\
&\leq E[g](u) - \liminf_{n \to +\infty} E[g - f](u_n) \\
&\leq E[g](u) - E[g - f](u) \\
&= E[f](u).
\end{aligned}$$

Which, with the help of (2.10), concludes the proof. $\square$

Here is a consequence of Theorem 2.13:

**Proposition 2.15.** *Let $f : \mathbb{R} \times \mathbb{R}^N \to [0, \infty)$ continuous, $p, q \in [1, \infty)$ such that $W^{1,p}(\Omega) \hookrightarrow L^q(\Omega)$, and let $\varphi \in W^{1,\infty}(\Omega)$. Assume that there exists $g : \mathbb{R} \times \mathbb{R}^N \to [0, \infty)$ globally continuous and convex with respect to the last variable, such that*

$$f(u, \xi) \leq g(u, \xi) \leq c(f(u, \xi) + |\xi|^p + |u|^q + 1), \qquad \forall u \in \mathbb{R}, \, \forall \xi \in \mathbb{R}^N,$$

*for some $c > 0$. Then for any $u \in W^{1,p}_\varphi(\Omega)$, there exists a sequence $(u_n) \subset W^{1,\infty}_\varphi(\Omega)$ such that*

$$u_n \to u \quad \text{strongly in} \quad W^{1,p}(\Omega),$$

*and*

$$\lim_{n \to +\infty} \int_\Omega f(u_n(x), \nabla u_n(x)) dx = \int_\Omega f(u(x), \nabla u(x)) dx.$$

*Proof.* In the case where $E[f](u) = +\infty$, this follows from the continuity of $f$ and the Fatou Lemma. Otherwise, we can apply Theorem 2.13. By assumption, $g$ is convex with respect to $\xi$ and furthermore

$$\int_\Omega g(u(x), \nabla u(x)) dx \leq c \left( \int_\Omega f(u(x), \nabla u(x)) dx + \|\nabla u\|_{L^p}^p + \|u\|_{L^q}^q + |\Omega| \right)$$
$$< +\infty,$$

which is what we needed. $\square$

Proposition 2.15 has a few interesting consequences:

**Corollary 2.16.** *Let $f : \mathbb{R} \times \mathbb{R}^N \to [0, \infty)$ continuous such that $f^{**}$ is continuous, $p, q \in [1, \infty)$ such that $W^{1,p}(\Omega) \hookrightarrow L^q(\Omega)$, and let $\varphi \in W^{1,\infty}(\Omega)$. Assume that*

$$f(u, \xi) \leq c(f^{**}(u, \xi) + |\xi|^p + |u|^q + 1), \qquad \forall u \in \mathbb{R}, \, \forall \xi \in \mathbb{R}^N,$$

*for some $c > 0$. Then for any $u \in W^{1,p}_\varphi(\Omega)$, there exists a sequence $(u_n) \subset W^{1,\infty}_\varphi(\Omega)$ such that*

$$u_n \to u \quad \text{strongly in} \quad W^{1,p}(\Omega),$$

*and*

$$\lim_{n \to +\infty} \int_\Omega f(u_n(x), \nabla u_n(x)) dx = \int_\Omega f(u(x), \nabla u(x)) dx.$$

*Proof.* It is sufficient to take $g := c(f^{**} + |\xi|^p + |u|^q + 1)$ in Proposition 2.15. $\square$

**Corollary 2.17.** *Let $f : \mathbb{R} \times \mathbb{R}^N \to [0, \infty)$ be a continuous function, $p, q \in [1, \infty)$ such that $W^{1,p}(\Omega) \hookrightarrow L^q(\Omega)$, and let $\varphi \in W^{1,\infty}(\Omega)$. Assume that $f$ is "convex at infinity"; that is, there exists $\tilde{f} : \mathbb{R} \times \mathbb{R}^N \to [0, \infty)$ convex in $\xi$ and continuous as well as $K > 0$ such that*
$$f(u, \xi) = \tilde{f}(u, \xi), \qquad \forall u \in \mathbb{R}, \, \forall \xi \in \mathbb{R}^N \setminus B_K.$$
*Assume furthermore that*
$$\sup\left\{ \left|f(u,\xi) - \tilde{f}(u,\xi)\right| \, : \, \xi \in B_K \right\} = O(|u|^q), \quad as \, |u| \to \infty. \tag{2.12}$$
*Then for any $u \in W^{1,p}_\varphi(\Omega)$, there exists a sequence $(u_n) \subset W^{1,\infty}_\varphi(\Omega)$ such that*
$$u_n \to u \quad strongly \, in \quad W^{1,p}(\Omega),$$
*and*
$$\lim_{n \to +\infty} \int_\Omega f(u_n(x), \nabla u_n(x)) dx = \int_\Omega f(u(x), \nabla u(x)) dx.$$

*Proof.* By (2.12) and continuity in $u$ of the term on its left-hand side, there exists $M > 0$ such that
$$\left|f(u,\xi) - \tilde{f}(u,\xi)\right| \leq M(|u|^q + 1), \qquad \forall u \in \mathbb{R}, \, \forall \xi \in \mathbb{R}^N.$$
Now let $g(u,\xi) := \tilde{f}(u,\xi) + M(|u|^q + 1)$ and apply Proposition 2.15. $\square$

The following result consists in the case where $f = f(\xi)$ is $C^2$ and not far from being convex (more specifically, the eigenvalues of the Hessian of $f$ are not too largely negative). Before stating it, here are some notations: given $f : \mathbb{R}^N \to [0, \infty)$ a $C^2$ function, we define the function $\lambda[f] : \mathbb{R}^N \to \mathbb{R}$ by
$$\lambda[f](\xi) := \inf_{|\zeta|=1} \langle \nabla^2 f(\xi) \cdot \zeta, \zeta \rangle, \qquad \forall \xi \in \mathbb{R}^N.$$

Hence $\lambda[f](\xi)$ is the smallest eigenvalues of the Hessian of $f$ at point $\xi$. Notice that $f$ being $C^2$, $\lambda[f]$ is continuous on $\mathbb{R}^N$. We can now state the result:

**Corollary 2.18.** *Let $f : \mathbb{R}^N \to [0, \infty)$ be a function of class $C^2$ and $p \in [1, \infty)$. Denote*
$$\theta(r) := \sup_{|\xi|=r} \left(\lambda[f](\xi)\right)^-, \qquad \forall r \geq 0. \tag{2.13}$$
*(Here we used the notation $a^- := \max(-a, 0)$). We assume that $\theta(r) = O(r^{p-2})$, as $r \to +\infty$. Moreover, if $p = 1$, we assume that $\int_0^{+\infty} \theta(r) dr < +\infty$. Then for any $u \in W^{1,p}_\varphi(\Omega)$, there exists a sequence $(u_n) \subset W^{1,\infty}_\varphi(\Omega)$ such that*
$$u_n \to u \quad strongly \, in \quad W^{1,p}(\Omega),$$
*and*
$$\lim_{n \to +\infty} \int_\Omega f(\nabla u_n(x)) dx = \int_\Omega f(\nabla u(x)) dx.$$

*Proof.* As always, we can assume that $\int_\Omega f(\nabla u)dx < +\infty$, otherwise the result follows by the Fatou Lemma. Notice first that if $\theta \equiv 0$, then $f$ is convex and Theorem 2.1 can immediately be applied. We assume in the following that $\theta$ is not identically 0 on $[0,\infty)$. Let $\gamma : [0,\infty) \to [0,\infty)$ be $C^2$ such that $\gamma(0) = \gamma'(0) = 0$ to be fixed later on, as well as $g_1 : \xi \mapsto \gamma(|\xi|)$ defined on $\mathbb{R}^N$. As such, $g_1$ is a $C^2$ radial function on $\mathbb{R}^N$. Let $\xi \in \mathbb{R}^N \setminus \{0\}$ and $\mathcal{B}$ an orthonormal basis of $\mathbb{R}^N$ having $\xi/|\xi|$ as its first vector, then a straightforward calculation using the radial property of $g_1$ gives that its Hessian quadratic form is represented in the basis $\mathcal{B}$ as the following diagonal matrix:

$$\mathrm{diag}\left(\gamma''(|\xi|),\ |\xi|^{-1}\gamma'(|\xi|),\ \ldots,\ |\xi|^{-1}\gamma'(|\xi|)\right),$$

and so,

$$\mathrm{spectrum}\left[\nabla^2 g_1(\xi)\right] = \{\gamma''(|\xi|),\ |\xi|^{-1}\gamma'(|\xi|)\}, \qquad \forall \xi \in \mathbb{R}^N \setminus \{0\}. \tag{2.14}$$

- If $p > 1$, choose $\gamma$ to be a $C^2$ function on $[0,\infty)$ such that $\gamma(0) = \gamma'(0) = 0$ and such that $\gamma(r) = r^p$ if $r \geq 1$. Then (2.14) gives the existence of some $c > 0$ such that it holds

$$\lambda[g_1](\xi) \geq c|\xi|^{p-2}, \qquad \forall \xi \in \mathbb{R}^N \setminus B_2.$$

To be clear, here $B_2$ refers to the centered ball of radius 2 in $\mathbb{R}^N$. Now using the assumption we made on $\theta$, one can find $R \geq 2$ and $M \geq 1$ such that,

$$\lambda[g_1](\xi) \geq M^{-1}\theta(|\xi|), \qquad \forall \xi \in \mathbb{R}^N \setminus B_R, \tag{2.15}$$

therefore, one gets

$$\lambda[f + Mg_1](\xi) \geq \lambda[f](\xi) + M\lambda[g_1](\xi) \geq 0, \qquad \forall \xi \in \mathbb{R}^N \setminus B_R. \tag{2.16}$$

If $f + Mg_1$ is convex on $\mathbb{R}^N$, then we can take $g := f + Mg_1$ and apply Theorem 2.13. Indeed, recall that $u \in W^{1,p}(\Omega)$ and thus that $\int_\Omega g_1(\nabla u)dx < +\infty$. Otherwise, because $f$ and $g_1$ are $C^2$ on $\mathbb{R}^N$, it holds from (2.16) that

$$0 > \inf_{\xi \in B_R} \lambda[f + Mg_1](\xi) =: -c > -\infty. \tag{2.17}$$

Let $g_2 : \mathbb{R}^N \to [0,\infty)$ be any $C^2$ convex function with linear growth satisfying that $\lambda[g_2] \geq c$, on $B_R$. Then, by (2.13), (2.15) and (2.17), $f + Mg_1 + g_2$ is convex on $\mathbb{R}^N$. Taking $g := f + Mg_1 + g_2$ in Theorem 2.13 gives the result.

- If $p = 1$, the proof is slightly different. Firstly, instead of taking $\gamma(r) = r$, we choose $\gamma$ such that $\gamma'' = \theta$ (as well as $\gamma(0) = \gamma'(0) = 0$). We claim that $\gamma''(r) = O(r^{-1}\gamma'(r))$, as $r \to +\infty$. Indeed, $\gamma'$ is a non decreasing function, converging to $\tilde{c} := \int_0^{+\infty} \theta > 0$. Thus $r^{-1}\gamma'(r) \sim \tilde{c}r^{-1}$. However, $\gamma''(r) = \theta(r) = O(r^{-1})$ by assumption. Let $R > 0$ and $M \geq 1$ be such that for every $r > R$, it holds that

$$\gamma''(r) \leq Mr^{-1}\gamma'(r). \tag{2.18}$$

Using the fact that $M \geq 1$, as well as (2.14) and (2.18), we obtain once again (2.15). Moreover, we have that $g_1$ has sublinear growth (because $\gamma'$ is bounded), thus again $\int_\Omega g_1(\nabla u)dx < +\infty$. From this point on, the proof follows just as in the case $p > 1$.

□

# 3 Generalization to the non-autonomous case

In this section we want to extend Theorem 2.4 to the non autonomous case, applying the results from [7]. Firstly we report [7, Theorem 3], the main conditions used in the paper are the following : we say that a measurable function $g : \Omega \times \mathbb{R} \times \mathbb{R}^N \to [0, \infty)$ satisfies condition $(\mathcal{H}_1)$ if for every $L = (L_1, L_2) \in (0, \infty)^2$, there exists a constant $C_L > 0$ such that for a.e. $x \in \Omega$ and every $(u, \xi) \in [-L_1, L_1] \times \mathbb{R}^N$, for every $\varepsilon > 0$, it holds

$$(g_\varepsilon^-)^{**}(x, u, \xi) \leq \frac{L_2}{\varepsilon^N} \quad \Rightarrow \quad g(x, u, \xi) \leq C_L(1 + (g_\varepsilon^-)^{**}(x, u, \xi)), \tag{$\mathcal{H}_1$}$$

where

$$g_\varepsilon^-(x, u, \xi) = \operatorname*{ess\,inf}_{y \in \Omega \cap B_\varepsilon(x)} g(y, u, \xi). \tag{3.1}$$

We also say that a measurable function $g : \Omega \times \mathbb{R} \times \mathbb{R}^N \to [0, \infty)$ satisfies condition $(\mathcal{H}_2)$ if there exists $\theta \in [1, \infty]$, $a \in L^\theta(\Omega)$ and a real number $u_0 > 0$ such that

$$g(x, u, 0) \leq a(x)|u|^{\frac{p^*}{\theta'}}, \qquad \text{for a.e. } x \in \Omega, \ \forall u \in \mathbb{R} \setminus [-u_0, u_0]. \tag{$\mathcal{H}_2$}$$

Here $p^*$ and $\theta'$ are respectively the Sobolev and Hölder conjugate exponents of $p$ and $\theta$. We now state the result proved in [7]:

**Theorem 3.1.** *Let $f : \Omega \times \mathbb{R} \times \mathbb{R}^N \to [0, \infty)$ be a Carathéodory function which is convex with respect to the last variable, and let $\varphi \in W^{1,\infty}(\Omega)$. Also assume that $f$ satisfies $(\mathcal{H}_1)$. If $N \geq 2$, we assume furthermore that $f$ satisfy $(\mathcal{H}_2)$. Then, for every $u \in W^{1,p}_\varphi(\Omega)$ such that $E[f](u) < +\infty$ there exists a sequence $(u_n) \subset W^{1,\infty}_\varphi(\Omega)$ such that $u_n \to u$ strongly in $W^{1,p}(\Omega)$ and*

$$\lim_{n \to +\infty} \int_\Omega f(x, u_n(x), \nabla u_n(x))dx = \int_\Omega f(x, u(x), \nabla u(x))dx.$$

**Hypothesis $\mathcal{C}$.** We will work with the following set of assumptions:

  a) $f$ satisfies Hypothesis $\mathcal{A}$,

  b) $f$ satisfies condition $(\mathcal{H}_1)$,

  c) if the dimension $N$ of $\Omega$ is strictly greater than 1, $f$ satisfies condition $(\mathcal{H}_2)$,

  d) $f^{**}$ is a Carathéodory function.

We introduce a technical lemma about an anti-jump condition to impose on the original Lagrangian $f$ to obtain the condition $(\mathcal{H}_1)$ for $f^{**}$.

**Lemma 3.2.** *Let $f : \Omega \times \mathbb{R} \times \mathbb{R}^N \to [0, \infty)$ be a measurable function. Then for any $\varepsilon > 0$,*

$$((f^{**})_\varepsilon^-)^{**} = (f_\varepsilon^-)^{**}.$$

The notation in this lemma is as introduced in (3.1)

*Proof.* Since $f^{**} \leq f$, then
$$((f^{**})_\varepsilon^-)^{**} \leq (f_\varepsilon^-)^{**} \tag{3.2}$$

To show the other inequality, we first recall that, using Proposition 1.9,
$$(f_\varepsilon^-)^{**}(x, u, \xi) = \inf\left\{\sum_i \alpha_i \left(\operatorname*{ess\,inf}_{y \in \Omega \cap B_\varepsilon(x)} f(y, u, \xi_i)\right) \;\Big|\; \sum_i \alpha_i \xi_i = \xi\right\},$$

and
$$(f^{**})_\varepsilon^-(x, u, \xi) = \operatorname*{ess\,inf}_{y \in \Omega \cap B_\varepsilon(x)} \left(\inf\left\{\sum_i \alpha_i f(y, u, \xi_i) \;\Big|\; \sum_i \alpha_i \xi_i = \xi\right\}\right).$$

For every convex combination we have that
$$\sum_i \alpha_i \operatorname*{ess\,inf}_{y \in \Omega \cap B_\varepsilon(x)} f(y, u, \xi_i) \leq \operatorname*{ess\,inf}_{y \in \Omega \cap B_\varepsilon(x)} \sum_i \alpha_i f(y, u, \xi_i),$$

and
$$\inf\left(\operatorname*{ess\,inf}_{y \in \Omega \cap B_\varepsilon(x)} \left\{\sum_i \alpha_i f(y, u, \xi_i) \;\Big|\; \sum_i \alpha_i \xi_i = \xi\right\}\right)$$
$$= \operatorname*{ess\,inf}_{y \in \Omega \cap B_\varepsilon(x)} \left(\inf\left\{\sum_i \alpha_i f(y, u, \xi_i) \;\Big|\; \sum_i \alpha_i \xi_i = \xi\right\}\right),$$

and so,
$$\inf\left\{\sum_i \alpha_i\left(\operatorname*{ess\,inf}_{y \in \Omega \cap B_\varepsilon(x)} f(y, u, \xi_i)\right) \;\Big|\; \sum_i \alpha_i \xi_i = \xi\right\}$$
$$\leq \operatorname*{ess\,inf}_{y \in \Omega \cap B_\varepsilon(x)} \left(\inf\left\{\sum_i \alpha_i f(y, u, \xi_i) \;\Big|\; \sum_i \alpha_i \xi_i = \xi\right\}\right),$$

that is, $(f_\varepsilon^-)^{**} \leq (f^{**})_\varepsilon^-$. Since $(f_\varepsilon^-)^{**}$ is a convex function with respect to $\xi$, we have
$$(f_\varepsilon^-)^{**} \leq ((f^{**})_\varepsilon^-)^{**},$$

thus with (3.2),
$$(f_\varepsilon^-)^{**} = ((f^{**})_\varepsilon^-)^{**}.$$

□

**Lemma 3.3.** *Let $f : \Omega \times \mathbb{R} \times \mathbb{R}^N \to [0, \infty)$ be a measurable function. Then if $f$ satisfies condition $(\mathcal{H}_1)$, $f^{**}$ satisfies it as well.*

*Proof.* By the previous lemma we can rewrite condition $(\mathcal{H}_1)$ for $f$ as
$$((f^{**})_\varepsilon^-)^{**}(x, u, \xi) \leq \frac{L_2}{\varepsilon^n} \Rightarrow f(x, u, \xi) \leq C_L(1 + ((f^{**})_\varepsilon^-)^{**}(x, u, \xi))$$

and since $f^{**} \leq f$ we have
$$((f^{**})_\varepsilon^-)^{**}(x, u, \xi) \leq \frac{L_2}{\varepsilon^n} \Rightarrow f^{**}(x, u, \xi) \leq C_L(1 + ((f^{**})_\varepsilon^-)^{**}(x, u, \xi)).$$

□

Now we can extend Theorem 2.4 to the non-autonomous case.

**Theorem 3.4.** *Let $f : \Omega \times \mathbb{R} \times \mathbb{R}^N \to [0, \infty)$ satisfy Hypothesis $\mathcal{C}$, and let $\varphi \in W^{1,\infty}(\Omega)$. Then for any $u \in W^{1,1}_\varphi(\Omega)$, there exists a sequence $(u_n) \subset W^{1,\infty}_\varphi(\Omega)$ such that*

$$u_n \to u \quad \text{in} \quad L^1(\Omega),$$

*and*

$$\lim_{n \to +\infty} \int_\Omega f(x, u_n(x), \nabla u_n(x)) dx = \int_\Omega f^{**}(x, u(x), \nabla u(x)) dx.$$

*Moreover,*

- *if there exists $\Phi : \mathbb{R}^N \to [0, \infty)$ superlinear such that*

$$f(x, u, \xi) \geq \Phi(\xi), \quad \text{for a.e. } x \in \Omega,\ \forall u \in \mathbb{R},\ \forall \xi \in \mathbb{R}^N,$$

  *then the sequence $(u_n)$ can be chosen so that $u_n \rightharpoonup u$ weakly in $W^{1,1}(\Omega)$;*

- *if for some $p \in (1, \infty)$ it holds that*

$$f(x, u, \xi) \geq c_1 |\xi|^p - c_2, \quad \text{for a.e. } x \in \Omega,\ \forall u \in \mathbb{R},\ \forall \xi \in \mathbb{R}^N,$$

  *for some $c_1, c_2 > 0$, and $u \in W^{1,p}(\Omega)$ then the sequence $(u_n)$ can be taken so that $u_n \rightharpoonup u$ weakly in $W^{1,p}(\Omega)$.*

*Proof.* Just as in the proof of Theorem 2.4, we can assume without loss of generality that $E[f^{**}](u) < +\infty$. By assumptions $f^{**}$ is a Carathéodory function and if $N \geq 2$, for $|u| \geq u_0$ and a.e. $x \in \Omega$,

$$f^{**}(x, u, 0) \leq f(x, u, 0) \leq a(x) |u|^{\frac{p^*}{\theta'}},$$

thus $f^{**}$ satisfies $(\mathcal{H}_2)$ in this case. By Lemma 3.3 and the fact that $f$ satisfies $(\mathcal{H}_1)$, for every $L = (L_1, L_2) \in (0, \infty)^2$, there exists a constant $C_L > 0$ such that for a.e. $x \in \Omega$ and every $(u, \xi) \in [-L_1, L_1] \times \mathbb{R}^N$, for every $\varepsilon > 0$ it holds

$$((f^{**})^-_\varepsilon)^{**}(x, u, \xi) \leq \frac{L_2}{\varepsilon^N} \Rightarrow f^{**}(x, u, \xi) \leq C_L (1 + ((f^{**})^-_\varepsilon)^{**}(x, u, \xi)).$$

Thus we can apply Theorem 3.1 to $f^{**}$ (which is a Carathéodory function by $\mathcal{C}$-d)) and so, for every $u \in W^{1,1}_\varphi(\Omega)$ there exists a sequence $(v_n) \subset W^{1,\infty}_\varphi(\Omega)$ such that $v_n \to u$ in $W^{1,1}(\Omega)$ and

$$\lim_{n \to +\infty} \int_\Omega f^{**}(x, v_n(x), \nabla v_n(x)) dx = \int_\Omega f^{**}(x, u(x), \nabla u(x)) dx.$$

By Theorem 1.19, for every $n$ there exists a sequence $(v^k_n) \subset W^{1,\infty}_\varphi(\Omega)$ such that $v^k_n \to u$ in $L^\infty(\Omega)$ as $k \to +\infty$, and

$$\lim_{k \to +\infty} \int_\Omega f(x, v^k_n(x), \nabla v^k_n(x)) dx = \int_\Omega f^{**}(x, u_n(x), \nabla u_n(x)) dx.$$

Thus, using a diagonal argument, there exists a sequence $(u_n) \subset W^{1,\infty}_\varphi(\Omega)$ such that $u_n \to u$ in $L^1(\Omega)$ and

$$\lim_{n \to +\infty} \int_\Omega f(x, u_n(x), \nabla u_n(x))dx = \int_\Omega f^{**}(x, u(x), \nabla u(x))dx.$$

If there exists $\Phi$ superlinear such that

$$\sup_n \int_\Omega \Phi(\nabla u_n(x))dx \leq \sup_n \int_\Omega f(x, u_n(x), \nabla u_n(x))dx < +\infty,$$

then, up to considering a subsequence, $(\nabla u_n)$ converges weakly in $L^1$ to some $v \in L^1$. It is easy to see (using an integration by parts argument) that actually $v = \nabla u$. So we have

$$u_n \rightharpoonup u \quad \text{weakly in} \quad W^{1,1}(\Omega).$$

If $p > 1$ the proof is the same taking $\Phi(\xi) = c_1|\xi|^p - c_2$. $\square$

At this point, we finally generalize Theorem 2.8 to the non-autonomous case.

**Theorem 3.5.** *Let $f : \Omega \times \mathbb{R} \times \mathbb{R}^N \to [0, \infty)$ satisfy Hypothesis C, and let $\varphi \in W^{1,\infty}(\Omega)$. Assume $u \in W^{1,p}_\varphi(\Omega)$ for some $p \in [1, \infty)$, and satisfies*

$$\int_\Omega f^{**}(x, u(x), \nabla u(x))dx = \int_\Omega f(x, u(x), \nabla u(x))dx.$$

*Then there exists a sequence $(u_n) \subset W^{1,\infty}_\varphi(\Omega)$ such that*

$$u_n \to u \quad \text{strongly in} \quad W^{1,p}(\Omega),$$

*and*

$$\lim_{n \to +\infty} \int_\Omega f(x, u_n(x), \nabla u_n(x))dx = \int_\Omega f(x, u(x), \nabla u(x))dx.$$

*Proof.* We work as in the proof of Theorem 2.8. For the case $p = 1$, we make use again of Lemma 2.7 to get the existence of some $\Phi : \mathbb{R}^N \to [0, \infty)$ superlinear and convex such that $\int_\Omega \Phi(\nabla u(x))dx < +\infty$. Then we define

$$g(x, u, \xi) := f(x, u, \xi) + \Phi(\xi) + \sqrt{1 + |\xi|^2}.$$

The only thing we have to show is that, under the assumption that $f$ satisfies $(\mathcal{H}_1)$ and $(\mathcal{H}_2)$, $g$ satisfies them as well. The remainder of the proof will follows then just as in Theorem 2.8. Notice that, by definition of $g_\varepsilon^-$ and the fact that $g - f$ does not depend on $x$ and $u$, it holds

$$g_\varepsilon^-(x, u, \xi) = f_\varepsilon^-(x, u, \xi) + \Phi(\xi) + \sqrt{1 + |\xi|^2}.$$

Therefore we have, from the fact that $\Phi$ was assumed to be convex, that

$$(g_\varepsilon^-)^{**}(x, u, \xi) \geq (f_\varepsilon^-)^{**}(x, u, \xi) + \Phi(\xi) + \sqrt{1 + |\xi|^2}, \tag{3.3}$$

by definition of $(g_\varepsilon^-)^{**}$. By Hypothesis $\mathcal{C}$, given $L = (L_1, L_2) \in (0, \infty)^2$, there exists a constant $C_L > 0$ such that for a.e. $x \in \Omega$ and every $(u, \xi) \in [-L_1, L_1] \times \mathbb{R}^N$, for every $\varepsilon > 0$, it holds

$$(f_\varepsilon^-)^{**}(x, u, \xi) \leq \frac{L_2}{\varepsilon^N} \quad \Rightarrow \quad f(x, u, \xi) \leq C_L(1 + (f_\varepsilon^-)^{**}(x, u, \xi)).$$

If for a.e. $x \in \Omega$ and every $(u, \xi) \in [-L_1, L_1] \times \mathbb{R}^N$ for every $\varepsilon > 0$,

$$(g_\varepsilon^-)^{**}(x, u, \xi) \leq \frac{L_2}{\varepsilon^N},$$

then $(f_\varepsilon^-)^{**}(x, u, \xi) \leq L_2/\varepsilon^N$ by (3.3), and thus

$$f(x, u, \xi) \leq C_L(1 + (f_\varepsilon^-)^{**}(x, u, \xi)).$$

So, taking $\tilde{C} = \max(1, C_L)$, we have by (3.3),

$$g(x, u, \xi) \leq C_L(1 + (f_\varepsilon^-)^{**}(x, u, \xi)) + \Phi(\xi) + \sqrt{1 + |\xi|^2}$$
$$\leq \tilde{C}(1 + (g_\varepsilon^-)^{**}(x, u, \xi)),$$

and thus $g$ satisfies $(\mathcal{H}_1)$. Now in the case $N \geq 2$, since $f$ satisfies condition $(\mathcal{H}_2)$ there exists $C' > 0$ such that

$$f(x, u, 0) + 1 + \Phi(0) \leq C'a(x)|u|^{\frac{p^*}{\theta'}}, \qquad \text{for a.e. } x \in \Omega, \ \forall u \in \mathbb{R} \setminus [-u_0, u_0].$$

Here, up to modify $a$ by taking $\tilde{a} := \max(a, 1)$, we assumed that the function $a$ in condition $(\mathcal{H}_2)$ satisfies $a(x) \geq 1$ on $\Omega$. By Theorem 1.16 we have that $g^{**}$ is a Carathéodory function since $g$ is superlinear with respect to the last variable. So $g$ satisfies Hypothesis $\mathcal{C}$, thus we can apply Theorem 3.4 to $g$ and argue as in the proof of Theorem 2.8. If $p > 1$ the proof is the same replacing $\Phi(\xi) + \sqrt{1 + |\xi|^2}$ by $|\xi|^p$.

□

In [7, Appendix A] there are some applications of Theorem 3.1 to the case $f(x, u, \cdot)$ nonconvex but dominated by a convex function $g$. The following result follows this idea:

**Theorem 3.6.** *Let $f : \Omega \times \mathbb{R} \times \mathbb{R}^N \to [0, \infty)$ be a Carathéodory function, and let $\varphi \in W^{1,\infty}(\Omega)$. Let $u \in W^{1,p}_\varphi(\Omega)$ for some $p \in [1, \infty)$ and assume that there exists $g : \Omega \times \mathbb{R} \times \mathbb{R}^N \to [0, \infty)$ a Carathéodory function satisfying the assumptions of Theorem 3.1 and $a \in L^1(\Omega, [0, \infty))$ such that*

$$f(x, u, \xi) \leq g(x, u, \xi) + a(x), \qquad \text{for a.e. } x \in \Omega, \ \forall u \in \mathbb{R}, \ \forall \xi \in \mathbb{R}^N,$$

*and*

$$\int_\Omega g(x, u(x), \nabla u(x))dx < +\infty.$$

*Then there exists a sequence $(u_n) \subset W^{1,\infty}_\varphi(\Omega)$ such that*

$$u_n \to u \quad \text{strongly in} \quad W^{1,p}(\Omega),$$

*and*

$$\lim_{n \to +\infty} \int_\Omega f(x, u_n(x), \nabla u_n(x))dx = \int_\Omega f(x, u(x), \nabla u(x))dx.$$

*Proof.* Applying Theorem 3.1 to $g$, there exists a sequence $(u_n) \subset W^{1,\infty}_\varphi(\Omega)$ converging to $u$ in $W^{1,p}(\Omega)$ such that

$$\lim_{n \to +\infty} \int_\Omega g(x, u_n(x), \nabla u_n(x)) + a(x) dx = \int_\Omega g(x, u(x), \nabla u(x)) + a(x) dx < +\infty.$$

The reminder of the proof follows just as in Theorem 2.14, using the fact that $f$ is Carathéodory and the Fatou Lemma. $\square$

# Acknowledgment

Firstly we want to thank our supervisors Professor Giulia Treu and Professor Pierre Bousquet for their important help in the writing of this paper. This research has been funded by the GNAMPA projects "Fenomeno di Lavrentiev, Bounded Slope Condition e regolarità per minimi di funzionali integrali con crescite non standard e lagrangiane non uniformemente convesse" (2024) and "Regolarità di soluzioni di equazioni paraboliche a crescita nonstandard degeneri" (2025), as well as from the University Research School EUR-MINT (State support managed by the National Research Agency for Future Investments program bearing the reference ANR-18-EURE-0023).